\documentclass[final,10pt]{elsarticle}
\usepackage{fullpage}

\usepackage{color}
\usepackage{multirow}
\usepackage{bm}

\usepackage{tikz}

\usetikzlibrary{arrows}
\tikzset{
  treenode/.style = {align=center, inner sep=2pt, text centered},
  arn_n/.style = {treenode, rectangle, rounded corners=1.5ex, black,  draw=black,
    fill=white, text width=8.5em},
  arn_r/.style = {treenode, rectangle, rounded corners=1.5ex, black,  draw=black,
    fill=white, text width=5.7em},
  arn_x/.style = {treenode, rectangle, rounded corners=1.5ex, black,  draw=black,
    fill=white, text width=6.8em},
}

\definecolor{Blue}{rgb}{0,0,1}
\definecolor{Red}{rgb}{1,0,0}

\newcommand{\spacebehavior}[1]{\ensuremath{\bm{y}_{#1}}} 
\newcommand{\spacebehaviormat}{\ensuremath{\bm{Y}}} 
\newcommand{\timebehavior}[1]{\ensuremath{\bm{z}_{#1}}} 
\newcommand{\timebehaviormat}{\ensuremath{\bm{Z}}} 
\newcommand{\noise}{\ensuremath{\bm{E}}} 
\newcommand{\timestepit}{\ensuremath{k}} 
\newcommand{\kmeans}{\ensuremath{\bar k}} 
\newcommand{\state}{\ensuremath{\bm{x}}} 
\newcommand{\statei}[1]{\ensuremath{x_{#1}}} 
\newcommand{\statesnap}{\ensuremath{\bm{X}}} 
 
\newcommand{\statesnaprowi}[1]{\ensuremath{\bm{x}^T_{#1}}} 
\newcommand{\statesnaptmp}{\ensuremath{\bar\statesnap}} 
\newcommand{\statesnapij}[2]{\ensuremath{{x}_{#1#2}}} 
\newcommand{\statesnaptmpij}[2]{\ensuremath{\bar x_{#1#2}}} 
\newcommand{\correlationmatrix}{\ensuremath{\bm{R}}} 
\newcommand{\f}{\ensuremath{\bm{f}}} 
\newcommand{\labelset}{\ensuremath{\bar E}} 
\newcommand{\ncluster}{\ensuremath{{n_c}}} 
\newcommand{\adjoint}{\ensuremath{\bm{y}}} 
\newcommand{\generic}{\ensuremath{\bm{w}}} 
\newcommand{\genericred}{\ensuremath{\hat\generic}} 
\newcommand{\genericredfine}{\ensuremath{\genericred^h}} 
\newcommand{\genericredcoarsetofine}{\ensuremath{\genericred_H^h}} 
\newcommand{\genericredcoarse}{\ensuremath{\genericred^H}} 
\newcommand{\residual}{\ensuremath{\bm{\tilde r}}} 
\newcommand{\A}{\ensuremath{\bm{A}}} 
\newcommand{\bvec}{\ensuremath{\bm{b}}} 
\newcommand{\nt}{\ensuremath{t}} 
\newcommand{\ntTrain}{\ensuremath{t_\mathrm{train}}} 
\newcommand{\outputs}{\ensuremath{z}} 
\newcommand{\outputfun}{\ensuremath{g}} 
\newcommand{\params}{\ensuremath{\bm{\mu}}} 
\newcommand{\param}[1]{\ensuremath{{\mu}_{#1}}} 
\newcommand{\leftbasis}{\ensuremath{\bm{W}}} 
\newcommand{\rightbasis}{\ensuremath{\bm{V}}} 
\newcommand{\rightbasisnew}{\rightbasis'} 
\newcommand{\barrightbasisfine}{\ensuremath{\bar\rightbasis^h}} 
\newcommand{\rightbasisfine}{\ensuremath{\rightbasis^h}} 
\newcommand{\rightbasisfineT}{\ensuremath{(\rightbasis^h)^T}} 
\newcommand{\rightbasiscoarse}{\ensuremath{\rightbasis^H}} 
\newcommand{\rightbasiscoarseT}{\ensuremath{\left(\rightbasis^H\right)^T}} 
\newcommand{\Q}{\ensuremath{\bm{Q}}} 
\newcommand{\R}{\ensuremath{\bm{R}}} 
\newcommand{\Qs}{\ensuremath{\bm{\bar Q}}} 
\newcommand{\Rs}{\ensuremath{\bm{\bar R}}} 
\newcommand{\Perm}{\ensuremath{\bm{\bar \Pi}}} 
\newcommand{\Permi}[1]{\ensuremath{\bm{\bar \pi}_{#1}}} 
\newcommand{\rightbasisi}[1]{\ensuremath{\bm{v}_{#1}}} 
 
\newcommand{\rightbasisfinei}[1]{\ensuremath{\rightbasisi{#1}^h}} 
\newcommand{\rightbasisfineiT}[1]{\left(\rightbasisfinei{#1}\right)^T} 
\newcommand{\rightbasisfineij}[2]{\ensuremath{v^h_{#1#2}}} 
\newcommand{\rightbasiscoarsei}[1]{\ensuremath{\rightbasisi{#1}^H}} 
\newcommand{\rightbasiscoarseij}[2]{\ensuremath{v^H_{#1#2}}} 
\newcommand{\nrb}{\ensuremath{{p}}} 
\newcommand{\canonical}[1]{\ensuremath{\bm{e}_{#1}}} 
\newcommand{\nrbavg}{\ensuremath{{\bar \nrb}}} 
\newcommand{\nrbinitial}{\ensuremath{{\nrb^{(0)}}}} 
\newcommand{\nrbfine}{\ensuremath{{q}}} 
\newcommand{\nominalrightbasis}{\ensuremath{\bm{V^{(0)}}}} 
\newcommand{\nominalrightbasisi}[1]{\ensuremath{\bm{v^{(0)}}_#1}} 
\newcommand{\nominalrightbasisij}[2]{\ensuremath{v^{(0)}_{#1#2}}} 
\newcommand{\redstate}{\ensuremath{\hat \state}} 
\newcommand{\redstatenew}{\ensuremath{\redstate'}} 
\newcommand{\redstatecoarse}{\ensuremath{\redstate^{H}}} 
\newcommand{\redstatecoarsek}{\ensuremath{(\redstate^{H})^\timestepit}} 
\newcommand{\redstatefine}{\ensuremath{\redstate^{h}}} 
\newcommand{\redstatefinek}{\ensuremath{(\redstate^{h})^\timestepit}} 
\newcommand{\redadjoint}{\ensuremath{\hat \adjoint}} 
\newcommand{\redadjointfine}{\ensuremath{\redadjoint^h}} 
\newcommand{\redadjointcoarsetofine}{\ensuremath{\redadjoint_H^h}} 
\newcommand{\redadjointcoarsetofineT}{\ensuremath{\left(\redadjoint_H^h\right)^T}} 
\newcommand{\redadjointcoarsetofinei}[1]{\ensuremath{\left[\hat
y_H^h\right]_{#1}}} 
\newcommand{\redadjointcoarse}{\ensuremath{\redadjoint^H}} 
\newcommand{\redadjointfineT}{\left(\redadjointfine\right)^T}

\newcommand{\errorIndicator}{\ensuremath{\bm{\delta}}} 
 
\newcommand{\errorIndicatorFine}{\ensuremath{\errorIndicator^h}} 
\newcommand{\errorIndicatorFinei}[1]{\ensuremath{\delta^h_{#1}}} 
 
\newcommand{\reference}{\ensuremath{\bar \state}} 
\newcommand{\basismap}{\ensuremath{\bm{A^n\left(\state;\params\right)}}} 
\newcommand{\modresidual}{\ensuremath{\bm{ r}}} 
\newcommand{\resetFreq}{\ensuremath{n_{\mathrm{reset}}}} 
\newcommand{\romtol}{\ensuremath{\epsilon_{\mathrm{ROM}}}} 
\newcommand{\fomtol}{\ensuremath{\epsilon}} 

\newcommand{\coarsetofine}{\ensuremath{\bm{I}_H^h}} 
 
\newcommand{\finetocoarse}{\ensuremath{\bm{I}_h^H}} 
\newcommand{\finetocoarseT}{\ensuremath{\left(\finetocoarseT\right)^T}} 

\newcommand{\childrenNo}{\ensuremath{{ C}}} 
 
\newcommand{\children}[1]{\ensuremath{{\childrenNo\left(#1\right)}}} 
 
\newcommand{\childreniNo}[1]{\ensuremath{{\childrenNo_{#1}}}} 
\newcommand{\childrenNoi}[1]{\ensuremath{{\childrenNo\left(#1\right)}}} 
\newcommand{\childrenNoij}[2]{\ensuremath{{\childrenNo\left(#1,#2\right)}}} 
\newcommand{\childrenij}[3]{\ensuremath{{\childrenNo_{#1}\left(#2,#3\right)}}} 
\newcommand{\nchildren}[1]{\ensuremath{{q_{#1}}}} 
\newcommand{\tmpchildreniNo}[1]{\ensuremath{{\bar \childrenNo_{#1}}}} 
\newcommand{\tmpchildreni}[2]{\ensuremath{{\bar\childrenNo_{#1}\left(#2\right)}}} 
\newcommand{\elementsNo}{\ensuremath{{ E}}} 
 
\newcommand{\elements}[1]{\ensuremath{\elementsNo\left(#1\right)}} 
\newcommand{\elementsiNo}[1]{\ensuremath{\elementsNo_{#1}}} 
\newcommand{\elementsi}[2]{\ensuremath{\elementsNo_{#1}\left(#2\right)}} 
 
\newcommand{\elementsNoij}[2]{\ensuremath{\elementsNo\left(#1,#2\right)}} 
\newcommand{\elementsNoi}[1]{\ensuremath{\elementsNo\left(#1\right)}} 
\newcommand{\iChild}{\ensuremath{j}} 
\newcommand{\tmpelementsi}[2]{\ensuremath{\bar\elementsNo_{#1}\left(#2\right)}} 
\newcommand{\tmpelementsiNo}[1]{\ensuremath{\bar\elementsNo_{#1}}} 
\newcommand{\nat}[1]{\ensuremath{\mathbb N(#1)}} 
\newcommand{\nodeindex}{\ensuremath{i}} 
\newcommand{\nodeindexa}{\ensuremath{j}} 
\newcommand{\nodeindexb}{\ensuremath{k}} 
\newcommand{\nnode}{\ensuremath{m}} 
\newcommand{\node}{\ensuremath{d}} 
\newcommand{\tmpnode}{\ensuremath{\bar d}} 
\newcommand{\ndof}{\ensuremath{n}} 
\newcommand{\nparams}{\ensuremath{n_{\params}}} 
\newcommand{\paramDomain}{\ensuremath{\mathcal D}} 
\newcommand{\paramtrain}{\ensuremath{\mathcal D_{\mathrm{train}}}} 
\newcommand{\powersetNo}{\ensuremath{\mathcal P}} 
\newcommand{\powerset}[1]{\ensuremath{\powersetNo\left(#1\right)}} 
\newcommand{\locerror}[2]{\ensuremath{\eta_{#1 #2}}} 
\newcommand{\flaggedvectors}{\ensuremath{ I}} 
\newcommand{\npartition}{\ensuremath{\alpha}} 
\newcommand{\basisweight}[1]{\ensuremath{\beta}_{#1}} 
\newcommand{\extravecCounter}{\ensuremath{l}} 
 
\newcommand{\extrabasis}[1]{\ensuremath{\bm{x}_{#1}}} 
 
\newcommand{\identity}[1]{\ensuremath{\bm{I}_{#1}}} 
\newcommand{\weightingMat}{\ensuremath{\bm{\Gamma}}} 
\newcommand{\weightingVec}[1]{\ensuremath{\bm{\gamma}_{#1}}} 
\newcommand{\rankbasis}{\ensuremath{r}} 
 
\newcommand{\urom}{\ensuremath{u_{\mathrm{ROM}}}}

\newcommand{\childgroupNo}{\ensuremath{ D}} 
\newcommand{\childgroup}[1]{\ensuremath{\childgroupNo_{#1}}} 
\newcommand{\childgroupvar}{\ensuremath{z}} 
\newcommand{\childgroupset}{\ensuremath{ K}} 

\newcommand{\vectorparenttochildNo}{\ensuremath{f}} 
\newcommand{\vectorparenttochild}[2]{\ensuremath{\vectorparenttochildNo\left(#1,#2\right)}} 
\newcommand{\vectorchildtoparentNo}{\ensuremath{f^{-1}}}

\newcommand{\nobs}{\ensuremath{{n_o}}} 
\newcommand{\burgersSize}{\ensuremath{{250}}} 
\newcommand{\burgersSizep}{\ensuremath{{251}}} 

\usepackage{mathtools}
\newcommand{\defeq}{\vcentcolon=}

\newcommand{\card}[1]{\ensuremath{\mathrm{card}\left(#1\right)}}
\newcommand{\range}[1]{\ensuremath{\mathrm{range}\left(#1\right)}}

\newcommand{\recentclusters}{\ensuremath{ D}} 
\newcommand{\recentclustersold}{\ensuremath{\bar \recentclusters}} 
\newcommand{\recentclustersoldi}[1]{\ensuremath{\recentclustersold\left(#1\right)}} 
\newcommand{\clustertosplit}{\ensuremath{d}} 
 
\newcommand{\clustercount}{\nnode} 
\newcommand{\icluster}{\ensuremath{i}} 
\newcommand{\kmeansmethod}{\ensuremath{\mathrm{kmeans}}} 
\newcommand{\timeVariable}{\ensuremath{\tau}} 
\newcommand{\tTrain}{\ensuremath{\tau_{\mathrm{train}}}}
\newcommand{\metric}{\ensuremath{\mathbf{\Theta}}}
\usepackage{algorithmic}

\newtheorem{theorem}{Theorem}
\newtheorem{lem}{Lemma}
\newenvironment{proof}[1][Proof]{\begin{trivlist}
\item[\hskip \labelsep {\bfseries #1}]}{\end{trivlist}}

\newcounter{lemctr}
\newcounter{pfctr}
\newcounter{remctr}
\newenvironment{remark}{\refstepcounter{remctr}\begin{trivlist}
\item {\emph {Remark}.\:}}
{$\blacksquare$\end{trivlist}}

\newcounter{examplectr}
\newenvironment{example}{\refstepcounter{examplectr}\begin{trivlist}
\item {\emph {Example}.\:}}
{$\blacksquare$\end{trivlist}}

\newcounter{propctr}
\newcounter{proposctr}
\newcommand{\RR}[1]{\ensuremath{\mathbb{R}^{ #1 }}}

\usepackage{subfigure, algorithmic, algorithm}
\usepackage{amsmath,amssymb}
\usepackage{graphics,epsfig}
\usepackage{wrapfigure}
\usepackage{subfigmat}
\usepackage{amsopn}

\begin{document}
\begin{frontmatter}
\title{Adaptive $h$-refinement for reduced-order models}
\author{Kevin Carlberg}
\ead{ktcarlb@sandia.gov}
\ead[url]{sandia.gov/~ktcarlb}
\address{Sandia National Laboratories}
\cortext[sandiacor]{7011 East Ave, MS 9159, Livermore, CA 94550. Sandia is a
multiprogram laboratory operated by Sandia
Corporation, a Lockheed Martin Company, for the United States Department of
Energy under contract DE-AC04-94-AL85000.}




\title{Adaptive $h$-refinement for reduced-order models}
\author{Kevin Carlberg}
	\address{Harry S.\ Truman Fellow, Quantitative Modeling
	\& Analysis Department\\
Sandia National Laboratories, P.O. Box 969, MS 9159, Livermore, CA
94551, USA.}


\begin{abstract}
This work presents a method to adaptively refine reduced-order models \emph{a
posteriori} without requiring additional full-order-model solves. The
technique is analogous to mesh-adaptive $h$-refinement: it enriches the
reduced-basis space online by `splitting' a given basis vector into several
vectors with disjoint support. The splitting scheme is defined by a tree
structure constructed offline via recursive $k$-means clustering of the state
variables using snapshot data. The method identifies the vectors to split
online using a dual-weighted-residual approach that aims to reduce error in an output
quantity of interest.  The resulting method generates a hierarchy of subspaces
online without requiring large-scale operations or full-order-model solves.
Further, it enables the reduced-order model to satisfy \emph{any prescribed
error tolerance} regardless of its original fidelity, as a completely refined
reduced-order model is mathematically equivalent to the original full-order
model.  Experiments on a parameterized inviscid Burgers equation highlight the
ability of the method to capture phenomena (e.g., moving shocks) not contained
in the span of the original reduced basis.

\end{abstract}

\end{frontmatter}
\section{Introduction}

Modeling and simulation of parameterized systems has become an essential tool across a wide range of
industries. However, the computational
cost of executing high-fidelity large-scale simulations is infeasibly high for
many time-critical applications. In particular, many-query
scenarios (e.g., sampling for solving statistical inverse problems) can require thousands of
simulations corresponding to different input-parameter instances of the
system; real-time contexts (e.g., model predictive control) require
simulations to execute in mere seconds.

Reduced-order models (ROMs) have been developed to mitigate this computational
bottleneck. First, they execute an `offline' stage during which
computationally expensive training tasks (e.g., evaluating the high-fidelity
model at several points in the input-parameter space) compute a representative
low-dimensional reduced basis for the system state. Then, during the
inexpensive `online' stage, these methods quickly compute approximate solutions for
arbitrary points in the input space via a projection process of the high-fidelity
full-order-model (FOM) equations onto the low-dimensional subspace spanned by
the reduced basis. They also introduce other approximations in the presence of
general (i.e., not low-order polynomial) nonlinearities.  See Ref.~\cite{surveyWillcoxGugercin} and references within
for a survey of current methods.

While reduced-order models almost always generate \emph{fast} online
predictions, there is no guarantee that they will generate sufficiently
\emph{accurate} online predictions.  In fact, the accuracy of online
predictions is predicated on the relevance of the training data to the online
problem: if a physical phenomenon was not observed during the offline stage, then
this feature will be missing from online predictions.  In general, the most
one can guarantee \emph{a priori} is that the ROM solution error is bounded by a prescribed
scalar over a finite set of `training points' in the input-parameter space
\cite{patera2007reduced}. While reduced-order models can be accurate at 
online points contained within a reasonable neighborhood of these training
points (see, e.g., Ref.~\cite{gunzburger2007reduced}), they are generally 
inaccurate for points far outside this set. 

This lack of error control\footnote{Note that reduced-order-model error bounds---which exist for
many problems---serve to quantify the error, while error control implies
reducing this error \emph{a posteriori}.} precludes ROMs from being employed in many
contexts.  For example, PDE-constrained optimization requires the
solution to
satisfy a prescribed forcing sequence to guarantee convergence
\cite{heinkenschloss2002analysis}. In uncertainty quantification, if the
epistemic uncertainty due to the ROM solution error dominates other sources of
uncertainty, the ROM cannot be exploited in a useful manner.  When
simulating parameterized highly nonlinear dynamical systems, it is unlikely
that any amount of training will fully encapsulate the range of complex
phenomena that can be encountered online; such problems require an efficient
refinement mechanism to generate accurate ROM predictions.

A few methods exist to improve a ROM solution when it is detected to be
inaccurate; however, they entail \emph{large-scale} operation counts.  The
most common approach is to revert to the high-fidelity model, solve the
associated high-dimensional equations for the current time step or
optimization iteration, add the solution to the reduced basis, and proceed
with the enriched reduced-order model
\cite{WeickumEldredJournal,arian2000trp,ryckelynck2005phm}. Another approach
adaptively improves the reduced-order model \emph{a posteriori} by generating
a Krylov subspace
\cite{carlberg-adaptive}; here, the reduced-order model serves to accelerate
the full-order solve to any specified tolerance. As our goal is to improve the
reduced-order model efficiently, i.e., without incurring large-scale
operations, none of these methods is appropriate.

Instead, this work proposes a novel approach inspired by mesh-adaptive
$h$-refinement. The main idea is to adaptively refine an inaccurate ROM
\emph{online} by
`splitting' selected reduced basis vectors into multiple vectors with disjoint
discrete support.  This splitting technique is defined by a tree structure
generated offline by
applying $k$-means clustering to the state variables. The method uses a
dual-weighted residual approach to select vectors to split online.
The resulting method generates a hierarchy of subspaces online
without requiring any large-scale operations or high-fidelity solves. Most
importantly, the methodology acts as a `failsafe' mechanism for the ROM:
$h$-adaptivity enables the ROM to satisfy \emph{any prescribed
error tolerance} online, as a fully refined ROM is mathematically equivalent to the original
full-order model under modest conditions. 

As a final note, some `adaptive' methods exist to tailor the ROM to specific
regions of the input space
\cite{Amsallem:2008bi,amsallem2009method,eftang2010hp,haasdonk2011training,drohmann2011adaptive,peherstorfer2013localized},
time domain \cite{drohmann2011adaptive,adaptiveTimePartitioning}, and state
space \cite{amsallemLocal,peherstorfer2013localized}.
However, these methods are primarily \emph{a priori} adaptive: they construct
separate ROMs for each region \emph{offline} with the goal of reducing the ROM
dimension.
While they can be used to improve the ROM \emph{a posteriori}, e.g., by
restarting the greedy algorithm online, doing so incurs additional
full-order-model solves, which is what we aim to avoid.

In the remainder of this paper,  matrices are denoted by capitalized bold
letters, vectors by lowercase bold letters, scalars by lowercase letters, and
sets by capitalized letters. The columns of a matrix $\bm{A}\in\RR{m\times
{\timestepit}}$ are denoted by $\bm{a}_i\in\RR{m}$, $i\in\nat{{\timestepit}}$
with $\nat{a}\defeq\{1,\ldots, a\}$ such that $\bm{A}\defeq\left[\bm{a}_1\
\cdots\ \bm{a}_{\timestepit}\right]$.  The scalar-valued matrix elements are
denoted by $a_{ij}\in\RR{}$ such that  $\bm{a}_j\defeq\left[a_{1j}\ \cdots\
a_{mj}\right]^T$, $j\in\nat{{\timestepit}}$.

\section{Problem formulation}
\subsection{Full-order model}
Consider solving a parameterized sequence of systems of equations
\begin{equation}\label{eq:FOMdyn}
\residual^{\timestepit}(\state^{{\timestepit}};\params) = 0
\end{equation}
for ${\timestepit}\in\nat{\nt}$, where
	 $\state^{\timestepit}\in\RR{\ndof}$ denotes the state at iteration
	${\timestepit}$, 
	$\params\in\paramDomain\subset\RR{\nparams}$ denotes the input parameters (e.g.,
	boundary conditions), $\residual^{\timestepit}:\RR{\ndof}\times\RR{\nparams}\rightarrow \RR{\ndof}$ denotes the
	residual operator at iteration ${\timestepit}$, and $\nt$ denotes maximum
	number of iterations. This formulation is quite general, as it
	describes, e.g., parameterized systems of linear equations ($\nt=1$,
	$\residual:(\state;\params)\mapsto \bvec(\params) - \A(\params)\state$) such
	as those arising from the finite-element discretization of elliptic PDEs, and
	parameterized ODEs $\dot \state
= \f(\state;\params)$ after time discretization by an implicit linear
multistep method (e.g., $\residual^{\timestepit}:(\state^{{\timestepit}};\params)\mapsto \state^{\timestepit} - \state^{{\timestepit}-1} - \Delta t
\f\left(\state^{\timestepit};\params\right)$
for the backward Euler scheme) such as those arising from the space- and
	time-discretization of parabolic and hyperbolic PDEs.
Assume that we are primarily interested in computing outputs
\begin{equation} \label{eq:outdyn}
\outputs^{\timestepit}  = \outputfun(\state^{\timestepit};\params)
\end{equation}
with $\outputs^{\timestepit}\in\RR{}$ and $\outputfun:\RR{\ndof}\times\RR{\nparams}\rightarrow\RR{}$.

When the dimension $\ndof$ is `large', computing the outputs of interest
$\outputs^\timestepit$ by first solving Eq.~\eqref{eq:FOMdyn} and subsequently
computing outputs via Eq.~\eqref{eq:outdyn} can be prohibitively expensive.
This is particularly true for many-query (e.g., statistical inversion) and
real-time (e.g., model-predictive control) problems that demand a fast
evaluation of the input--output map
$\params\mapsto\{\outputs^{1},\ldots,\outputs^{\nt}\}$. 

\subsection{Reduced-order model}
Model-reduction techniques aim to reduce the burden of solving Eq.\
\eqref{eq:FOMdyn} by employing a projection process. 
First, they execute a computationally expensive offline stage 
(e.g., solving Eq.~\eqref{eq:FOMdyn} for a training
set $\params\in\paramtrain\subset\paramDomain$) to construct 1) a
low-dimensional trial basis (in matrix form) $\rightbasis\in
\RR{\ndof\times\nrb}$ with $\nrb\ll\ndof$ that (hopefully) captures the behavior of
the state $\state$ throughout the parameter domain $\paramDomain$, and 2) an associated test basis $\leftbasis\in\RR{\ndof\times\nrb}$.
Then, during the computationally inexpensive online stage, these methods
approximately solve Eq.~\eqref{eq:outdyn} for arbitrary $\params\in
\paramDomain$ by 
searching for solutions
in the trial subspace $\reference + \range{\rightbasis}\subset\RR{\ndof}$ (with
$\reference\in\RR{\ndof}$ a chosen reference configuration) and
enforcing the residual $\residual^{\timestepit}$ to be orthogonal to the test subspace
$\range{\leftbasis}\subset\RR{\ndof}$:
\begin{equation}\label{eq:ROMdyn}
\leftbasis^T\residual^{\timestepit}(\reference +
\rightbasis\redstate^{{\timestepit}};\params) = 0.
\end{equation}
Here, $\redstate^{\timestepit}\in\RR{\nrb}$ denotes the generalized
coordinates of the reduced-order-model solution $\reference +
\rightbasis\redstate^{\timestepit}$ at iteration $\timestepit$.
When the residual operator exhibits general nonlinear dependence on the state
or is non-affine in the
inputs, additional
complexity-reduction approximations such as empirical interpolation
\cite{barrault2004eim}, collocation
\cite{LeGresleyThesis,astrid2007mpe,ryckelynck2005phm}, discrete
empirical interpolation
\cite{chaturantabut2010journal,drohmannEOI},  or gappy proper orthogonal
decomposition (POD) \cite{astrid2007mpe,carlbergJCP} are required to ensure that computing the low-dimensional
residual $\leftbasis^T\residual^{\timestepit}$ incurs an $\ndof$-independent
operation count. For simplicity, we do not consider such approximations in
the present work; future work will entail extending the proposed method to
such `hyper-reduced' order models.

In many cases, the test basis can be expressed as $\leftbasis = \basismap\rightbasis$. For example, $\basismap = I$ for
Galerkin projection; balanced
truncation uses $\basismap = \bm{Q}$, where $\bm{Q}$ is the observability
Gramian of the linear time-invariant system; the least-squares
Petrov--Galerkin projection \cite{LeGresleyThesis,carlbergJCP}
underlying the GNAT method employs $\basismap = {\partial
\residual^{\timestepit}}/{\partial \state}\left(\state,\params\right)$;
for linearized compressible-flow problems, 
$\basismap$ can be chosen to guarantee stability \cite{barone2009stable}. When
this holds, the
Petrov--Galerkin projection \eqref{eq:ROMdyn} is equivalent to a Galerkin
projection performed on the modified residual $\modresidual^{\timestepit} \defeq
\basismap^T\residual^{\timestepit}$:
\begin{equation}\label{eq:ROMdyn2}
\rightbasis^T\modresidual^{\timestepit}(\reference + \rightbasis\redstate^{{\timestepit}};\params) =
0,
\end{equation}
for $\timestepit\in\nat{\nrb}$. In the remainder of this paper, Eq.~\eqref{eq:ROMdyn2} will be considered the
governing equations for the reduced-order model.

\subsection{Objective: adaptive refinement}
The goal of this work is as follows: given a reduced basis $\rightbasis$ and
online ROM solution $\redstate^{\timestepit}$ to Eq.~\eqref{eq:ROMdyn2} for
iteration $\timestepit$, 1) detect 
if the solution is sufficiently accurate, 2) if it is not sufficiently
	accurate, efficiently generate a
higher-dimensional reduced basis $\rightbasisnew$ with $\range{\rightbasis}\subseteq\range{\rightbasisnew}$ in a goal-oriented manner
that aims to reduce errors
in the output $\outputs^{\timestepit}$, 3) compute an associated
solution $\redstatenew^{\timestepit}$, 4) repeat until desired accuracy is reached.

To generate this hierarchy of subspaces efficiently, we propose an analogue to
adaptive $h$-refinement, wherein selected basis vectors $\rightbasisi{i}$ are
`split' online into multiple vectors with disjoint support (i.e., the element
set with nonzero entries). Like all $h$-refinement techniques, the proposed
method consists of the following components:
\begin{enumerate} 
\item \emph{Refinement mechanism}. In typical $h$-refinement, this is defined
by the mesh-refinement method applied to finite elements or volumes. The proposed
method refines the solution space by splitting the support of the basis
vectors using a tree structure constructed via $k$-means clustering of the
state variables. Section \ref{sec:refinement} describes this component.
\item \emph{Error indicators}. Goal-oriented methods for $h$-refinement often
1) solve a coarse dual problem, 2) prolongate the adjoint solution to a
representation on the fine grid, and 3) compute error estimates of the output using
first-order analysis. The proposed method employs an analogous goal-oriented
dual-weighted residual approach. Section \ref{sec:errorInd} presents this.
\item \emph{An adaptive algorithm}. The proposed algorithm identifies when
refinement is required online and employs error indicators decide on the particular
refinement, i.e., which basis vectors should be refined, and how they should
be refined. Section \ref{sec:adaptive} provides this algorithm.
\end{enumerate}

\section{Refinement mechanism}\label{sec:refinement}
The method assumes that an initial reduced basis
$\nominalrightbasis\in\RR{\ndof\times \nrbinitial}$ is provided, which is
subsequently `split' to add fidelity to the ROM online. Section \ref{sec:tree} describes the tree data
structure that constitutes the splitting mechanism, Section \ref{sec:basisSplit} describes
how this mechanism leads to an algebraic refinement strategy, Section \ref{sec:properties}
highlight critical properties of the refinement method, and Section \ref{sec:kmeans} describes
construction of the tree via $k$-means clustering.

\subsection{Tree data structure}\label{sec:tree}
To begin, we define a tree data structure that characterizes the refinement
mechanism. The tree is characterized by a \emph{child} function
$\childrenNo:\nat{\nnode}\rightarrow\powerset{\nat{\nnode}}$ that describes
the topology of the tree and an
\emph{element} function
$\elementsNo:\nat{\nnode}\rightarrow\powerset{\nat{\ndof}}$ that describes the
set of nonzero vector entries associated with each tree node. Here,
$\nnode$ denotes the number of nodes in the tree and $\powersetNo$
denotes the powerset.

Each basis vector $\rightbasisi{i}$, $i\in\nat{\nrb}$ is characterized by a 
particular node on the tree $\node_i\in\nat{\nnode}$, a set of nonzero entries (i.e.,
support) 
$\elementsNoi{\node_i}$, and possible splits 
$\childrenNoi{\node_i}$.
If a given vector $\rightbasisi{i}$ is split, then it
is replaced in the basis by $\nchildren{i}\defeq \card{\children{\node_i}}$
child vectors whose 
set of nonzero entries is defined by $\elementsNoi{k}$, $k\in\children{\node_i}$; the
values of these nonzero entries are the same as those of the original vector
$\rightbasisi{i}$.


We enforce the following conditions for the tree:
\begin{enumerate} 
\item\label{req:globalSupport} The root node includes all elements: $\elements{1}
= \nat{\ndof}$, which 
is consistent with the possibly global support of the original reduced basis
$\nominalrightbasis$.
\item\label{req:childrenParent} The children have disjoint support, and the union of their
support equals that of the parent: For all $i\in\nat{\nnode}$, 
\begin{gather} 
\label{eq:reqDisjoint}\elements{\nodeindexa}\cap\elements{\nodeindexb}=\emptyset,\quad\forall
\nodeindexa,\nodeindexb\in\children{\nodeindex},\ \nodeindexa\neq\nodeindexb
\\
\label{eq:reqChildren}\mathop{\cup}\limits_{\nodeindexa\in\children{\nodeindex}}\elements{\nodeindexa}=\elements{\nodeindex}.
\end{gather} 
\item\label{req:fullSplit} Each element is associated with a single leaf node:
\begin{equation} \label{eq:requireall}
\forall l\in\nat{\ndof},\ \exists \nodeindex\in\nat{\nnode}\ |\
\elements{\nodeindex}
= l,\ \children{\nodeindex} = \emptyset.
\end{equation} 
\end{enumerate}
As will be shown, these requirements guarantee several critical properties of
the method.

\begin{example}
Consider an example with $\ndof=6$ and an initial reduced basis
$\nominalrightbasis=\nominalrightbasisi{1}$ of dimension
1.  Figure \ref{fig:treeEx} depicts an example of a tree structure
for this case.

 \begin{figure}[htbp] 
  \centering 
\begin{tikzpicture}[->,>=stealth', level 1/.style={sibling distance=6.5cm}
,
level 2/.style={sibling distance=2.5cm},
level 3/.style={sibling distance=2.3cm}]
\node [arn_n] {$\node=1$\\ $\children{1} = \{2,3\}$\\$\elements{1} = \{1,\ldots,6\}$}
    child{ node [arn_n] {$\node=2$\\ $\children{2} = \{4,5,6\}$\\$\elements{2} =
		\{1,3,4\}$} 
            child{ node [arn_r] {$\node=4$\\ $\children{4} = \emptyset$\\$\elements{4} = \{1\}$} 
            }
            child{ node [arn_r] {$\node=5$\\ $\children{5} = \emptyset$\\$\elements{5} = \{3\}$}
            }                            
            child{ node [arn_r] {$\node=6$\\ $\children{6} = \emptyset$\\$\elements{6} = \{4\}$}
            }                            
    }
    child{ node [arn_n] {$\node=3$\\ $\children{3} = \{7,8\}$\\$\elements{3} = \{2,5,6\}$}
            child{ node [arn_r] {$\node=7$\\ $\children{7} = \emptyset$\\$\elements{7} = \{2\}$}
            }                            
            child{ node [arn_x] {$\node=8$\\ $\children{8} = \{9,10\}$\\$\elements{8} = \{5,6\}$}
							child{ node [arn_r] {$\node=9$\\ $\children{9} = \emptyset$\\$\elements{9} = \{5\}$}
							}
							child{ node [arn_r] {$\node=10$\\ $\children{10} = \emptyset$\\$\elements{10} = \{6\}$}
							}
            }                            
		}
; 
\end{tikzpicture}
	  \caption{Tree example with $\ndof=6$} 
		 \label{fig:treeEx} 
		  \end{figure}
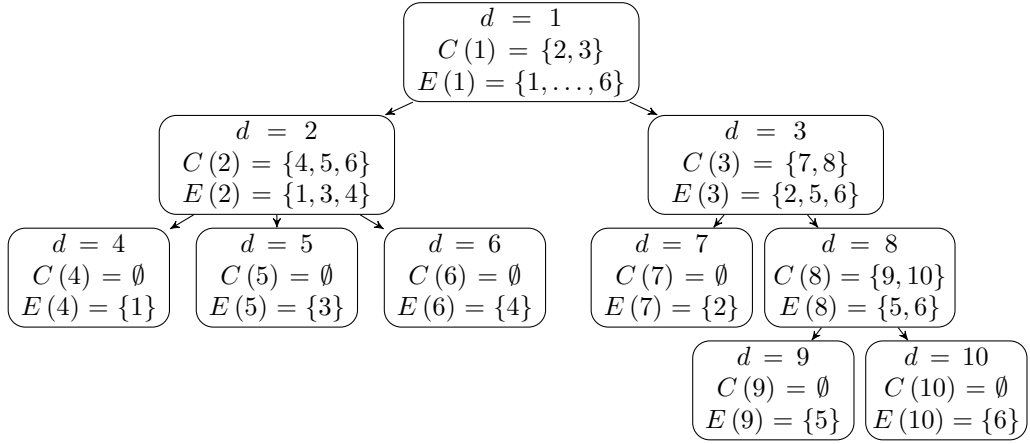 
Suppose the basis has been split into $\nrb=4$ according to the tree in Figure
\ref{fig:treeEx} with $\node_1 = 2$, $\node_2 = 7$, $\node_3 = 9$, and
$\node_4=10$; then, the refined reduced basis is
\begin{equation} 
\rightbasis = \left[\begin{array}{c c c c}
\nominalrightbasisij{1}{1} & 0 &0 & 0\\
0 & \nominalrightbasisij{2}{1} &0 & 0\\
\nominalrightbasisij{3}{1} & 0 &0 & 0\\
\nominalrightbasisij{4}{1} & 0 &0 & 0\\
0 & 0 &\nominalrightbasisij{5}{1} & 0\\
0 & 0 &0&\nominalrightbasisij{6}{1} \\
\end{array}
\right].
\end{equation} 
\end{example}

In the sequel, we overload the child function for the two-argument case such that
	$\childrenNoij{i}{j}$ denotes the $j$th child node of parent node $i$,
	where ordering of the children is implied by the binary relation $\leq$ on
	the natural numbers. Similarly, the overloaded element function
	$\elementsNoij{i}{j}$ is the $j$th element for node $i$; again,
	ordering of the elements is implied by the relation $\leq$ on the
	natural numbers.

\subsection{Refinement via basis splitting}\label{sec:basisSplit}
We now put the basis-splitting methodology in the framework of typical
$h$-refinement techniques.
First, define a `coarse' basis $\rightbasiscoarse\in\RR{\ndof\times\nrb}$,
which
is initially equal to the
nominal basis $\nominalrightbasis\in\RR{\ndof\times\nrbinitial}$ with
$\nrbinitial\leq\nrb$. As this initial basis may have global support, it is characterized by
$\node_i = 1$, $i\in\nat{\nrbinitial}$; this is permissible due to Condition
\ref{req:globalSupport} of Section \ref{sec:tree}.
Also define a `fine' basis
corresponding to the coarse basis with all vectors split 
according to the
children of the current node.
We can express the relationship between the coarse and fine
bases as
 \begin{equation} \label{eq:rightbasisfineSimple}
 \rightbasiscoarse = \rightbasisfine \coarsetofine,
  \end{equation} 
	where $\rightbasisfine\in\RR{\ndof\times \nrbfine}$ with $\nrbfine\geq\nrb$
	denotes the fine basis and $\coarsetofine\in\{0,1\}^{\nrbfine \times\nrb}$ denotes the
	\emph{prolongation} operator. 
Then, for any generalized coordinates $\genericredcoarse\in\RR{\nrb}$
associated with the coarse basis
$\rightbasiscoarse$, we can compute the corresponding fine
representation $\genericredfine\in\RR{\nrbfine}$ associated with the fine
basis $\rightbasisfine$ as 
\begin{equation}\genericredcoarsetofine = \coarsetofine\genericredcoarse,
\end{equation}
which ensures that $\rightbasiscoarse\genericredcoarse =
\rightbasisfine\genericredcoarsetofine$. 
Note this prolongation operator is exact, unlike typical mesh-refinement
strategies, where this operator is often defined as a linear or quadratic
interpolant of the coarse solution on the fine grid.
The \emph{restriction} operator is not
uniquely defined, but can be set, e.g., to
\begin{equation}\finetocoarse = \left(\coarsetofine\right)^+,
\end{equation}
where the superscript $+$ denotes the Moore--Penrose pseudoinverse.

 Using the tree structure defined in Section \ref{sec:tree}, we can precisely
 define these quantities. 
We first introduce
	 the mapping $\vectorparenttochildNo:(i,j)\mapsto k$, which provides the fine
	basis-vector index $k$ corresponding to the $j$th child of the $i$th coarse basis
	vector. We define it as
\begin{equation} 
\vectorparenttochild{i}{j} = 
\sum\limits_{k <
 i}\nchildren{k} +
 j, \quad j\in\nat{\nchildren{i}},\quad i \in\nat{\nrb}.
 \end{equation} 
 In particular, note that if node $\node_i$  is a leaf (i.e.,
 $\children{\node_i} = \emptyset$), then $\vectorparenttochild{i}{j}$ does not exist for any
 $j$.
 Similarly, the inverse mapping $\vectorchildtoparentNo:k\mapsto(i,j)$ 
 yields the coarse basis-vector index $i$ and child index $j$ corresponding to
 fine basis vector $k$.

 Now, the number of vectors in the fine reduced basis
 is simply 
\begin{equation} 
 \nrbfine = \sum_{i=1}^\nrb\nchildren{i}.
 \end{equation}
From Condition \ref{req:childrenParent} of Section \ref{sec:tree}, we can
write the fine reduced basis as
 \begin{align} \label{eq:rightbasisfine}
 \rightbasisfineij{i}{j}=
 \begin{cases}
 \rightbasiscoarseij{i}{l},\quad &\exists k\ | \ j = \vectorparenttochild{l}{k}, \
 i\in\elements{\childrenNoij{\node_l}{k}}\\
0,\quad &\mathrm{otherwise}.
 \end{cases}
  \end{align} 
	 and the prolongation operator induced by the proposed splitting scheme as 
 \begin{equation} \label{eq:coarsetofine}
 [\coarsetofine]_{ij} = 
 \begin{cases}
 1,\quad \exists k\ | \ i = \vectorparenttochild{j}{k}\\
0,\quad \mathrm{otherwise}.
 \end{cases}
  \end{equation} 

\subsection{Properties}\label{sec:properties}
This section highlights several key properties of this refinement method.
\begin{lem}[Hierarchical subspaces]\label{lem:hierarchy}
The method generates a hierarchy of subspaces such that
	$\range{\rightbasiscoarse}\subseteq\range{\rightbasisfine}$.
	\begin{proof}
	This result is self-evident from Eq.~\eqref{eq:rightbasisfineSimple}, as
 \begin{equation} 
 \range{\rightbasiscoarse} = \{\rightbasisfine \generic\ |\ \generic
 \in\range{\coarsetofine}\subseteq\RR{\nrbfine}\}\subseteq\{\rightbasisfine \generic\ |\ \generic
 \in\RR{\nrbfine}\} = \range{\rightbasisfine}.
  \end{equation} 
	\end{proof}
\end{lem}
\begin{theorem}[Monotonic convergence]\label{thm:monotone}
If the reduced-order model
\eqref{eq:ROMdyn2} is \textit{a priori} convergent, i.e., its solution satisfies
\begin{equation} \label{eq:optProbROM}
\rightbasis\redstate^{\timestepit} =
\arg\min_{\generic\in\range{\rightbasis}}\|\state^{\timestepit} -\reference-
\generic\|_{\metric},
\end{equation} 
for some norm $\|\cdot\|_\metric$,
then
the proposed refinement method guarantees monotonic convergence of the
reduced-order-model solution, i.e.,
\begin{equation} 
\|\state^{\timestepit} -\reference-
\rightbasisfine\redstatefinek\|_{\metric}
\leq
\|\state^{\timestepit} -\reference-
\rightbasiscoarse\redstatecoarsek\|_{\metric}.
\end{equation} 
\begin{proof}
This follows directly from Lemma \ref{lem:hierarchy}, as the
coarse-basis solution is contained in the span of the fine basis
$\rightbasiscoarse\redstatecoarsek\in\range{\rightbasiscoarse}\subseteq\range{\rightbasisfine}$.
\end{proof}
\end{theorem}
One example of a reduced-order model that satisfies the conditions of Theorem
\ref{thm:monotone} arises when the residual is linear in the
state and its Jacobian 
${\partial
\modresidual^{\timestepit}}/{\partial \state}\left(\params\right)$ is symmetric and positive
definite. In this case, a Galerkin-projection ROM satisfies Eq.~\eqref{eq:optProbROM}
for $\metric = {\partial
\modresidual^{\timestepit}}/{\partial \state}\left(\params\right)$ with
$\|\generic\|_{{\partial
\modresidual^{\timestepit}}/{\partial \state}\left(\params\right)} \defeq \sqrt{\generic^T{\partial
\modresidual^{\timestepit}}/{\partial \state}\left(\params\right)\generic}$.
Another example is
least-squares
Petrov--Galerkin applied to a parametrized system of linear equations
\cite{carlbergJCP}, where $\metric = \left({\partial
\modresidual^{\timestepit}}/{\partial \state}\left(\params\right)\right)^T{\partial
\modresidual^{\timestepit}}/{\partial \state}\left(\params\right)$.
\begin{theorem}[Convergence to the full-order model]\label{thm:fomconv}
If every element has a nonzero entry in one of the original
reduced-basis vectors, i.e.,
\begin{equation} \label{eq:globalSupport}
\forall l\in\nat{\ndof},\ \exists\ (i,j)\in\nat{\ndof}\times\nat\nrbinitial\ |\
\nominalrightbasisij{i}{j}\neq 0,
\end{equation} 
and Eq.~\eqref{eq:requireall} holds, then a completely split basis
yields a reduced-order
model equivalent to the full-order model.
\begin{proof}
Under these conditions, a completely split basis can be written as
$\rightbasis\in\RR{\ndof\times{\ndof\nrbinitial}}$ with all basis vectors in
the leaf-node state, i.e., 
$\childrenNoi{\node_i} = \emptyset$, $i\in\nat{\ndof\nrbinitial}$. Because
Eq.~\eqref{eq:requireall} guarantees that each element is associated with a
single leaf node, this implies that 
\begin{equation} \label{eq:globalSupport}
\forall l\in\nat{\ndof},\ \exists\ i\in\nat{\ndof\nrbinitial}\ |\
\rightbasisi{i}=\canonical{l}\basisweight{i},
\end{equation} 
where $\canonical{l}\in\{0,1\}^\ndof$ denotes the $l$th canonical unit
vector and $\basisweight{i}\neq 0$, $i\in\nat{\ndof\nrbinitial}$.
Eq.~\eqref{eq:globalSupport} implies that the completely split basis can be
post-multiplied by a (weighted) permutation matrix to yield the
$\ndof\times\ndof$ identity matrix $\identity{\ndof}$, i.e.,
\begin{equation} \label{eq:splitBasisToId}
\identity{\ndof} = \rightbasis\weightingMat.
\end{equation} 
Here, the matrix $\weightingMat\in\RR{\ndof\nrbinitial\times \ndof}$ consists
of columns
 \begin{equation} 
 \weightingVec{l} = \frac{1}{\basisweight{i}}\canonical{i},\quad i\in\{j\ |\
 \rightbasisi{j}= \canonical{l}\basisweight{j}\},\quad l\in\nat{\ndof}.
  \end{equation} 
Eq.~\eqref{eq:splitBasisToId} implies that 
 \begin{equation} 
 \range{\identity{\ndof}} =
 \RR{\ndof}\subseteq\range{\rightbasis}\subseteq\RR{\ndof},
  \end{equation} 
	which completes the proof.
\end{proof}
\end{theorem}
Lemma \ref{lem:hierarchy} and Theorem \ref{thm:fomconv} show that the proposed
refinement method enables the reduced-order model to generate a sequence of
hierarchical subspaces that converges to the full-order model under modest
assumptions. Thus, the method acts as a `failsafe' mechanism: it allows the
reduced-order model to generate \emph{arbitrarily accurate solutions}.
Despite this result, the associated rate of convergence is unknown, which
precludes any \emph{a priori} guarantee that the $h$-adaptive ROM 
will remain truly low dimensional for stringent accuracy
requirements. However, numerical experiments in Section \ref{sec:experiments} demonstrate that 
the proposed method often leads to accurate responses with
low-dimensional refined bases.

\begin{remark}
Note that the refinement method does not preclude a rank-deficient basis; this
can be seen from Theorem \ref{thm:fomconv}, wherein a completely split basis
has $\ndof\nrbinitial\geq \ndof$ columns. To detect (and remove) rank deficiency, the refinement
algorithm computes a rank-revealing QR factorization after each split  (Steps
\ref{step:qr}--\ref{step:fullrankstart} of Algorithm \ref{alg:refine} and Steps
\ref{step:qr2}--\ref{step:fullrankstart2} of Algorithm \ref{alg:refinemult}).
\end{remark}

\subsection{Tree construction via $k$-means clustering of the state variables}\label{sec:kmeans}
Any tree that satisfies Conditions
\ref{req:globalSupport}--\ref{req:fullSplit} of Section \ref{sec:tree} will
lead to the critical properties proved in Section \ref{sec:properties}. This
section presents one such tree-construction approach, which executes offline and employs the
following heuristic: 
\begin{quote}
State variables $\statei{i}$ that tend to be
strongly positively or negatively correlated can be accurately represented by
the same generalized coordinate, and should therefore reside in the same tree
node. 
\end{quote}
\begin{example}
To justify this heuristic, consider an example with $\ndof = 6$ degrees of freedom and $\nobs =
8$ observations of the state, e.g., from a computed time history. Assume
that snapshots can be decomposed as
 \begin{equation} 
 \statesnap = \sum_{i=1}^3\spacebehavior{i}\timebehavior{i}^T + 0.1\noise
  \end{equation} 
	where $\noise\in\left[-1,1\right]^{\ndof\times \nobs}$ is a matrix of
random uniformly distributed noise and the data matrices are
 \begin{equation*} 
 \timebehaviormat = \left[\begin{array}{c c c c c c c c c c}
   -2.2083  & -5.1072   & 2.6816&    9.3277&   -6.4506&   -3.2548&    4.2237&   -3.2557\\
   -2.9810  &  0.6557   & 3.0474&    5.5252&    2.7674&    2.3311&    9.6190&   -6.6484\\
   -2.4547  &  5.2676   &-3.6434&    5.5661&   -7.5449&    9.3079&   -2.0459&   -0.0728\\
 \end{array}
\right]^T
  \end{equation*} 
 
 \begin{equation*} 
 \spacebehaviormat = \left[\begin{array}{c c c c c c c c}
   -3.9885 &        0&         0&         0&         0&         0\\
         0 &        0&    8.6843&         0&         0&   -1.6393\\
         0 &  -1.7288&         0&    6.0559&    2.2407&         0\\
 \end{array}
\right]^T.
  \end{equation*} 
The sparsity structure of $\spacebehaviormat$ implies that the following sets of
state variables are strongly correlated or anti-correlated across observations: $\{1\}$,
$\{3,6\}$, and $\{2,4,5\}$. This is apparent from computing the matrix of
sample correlation coefficients:
\begin{equation}
\correlationmatrix =\left[
 \begin{array}{c c c c c c}
    1.0000 &   0.1526&   -0.5698&   -0.1534&   -0.1554&    0.5705\\
    0.1526 &   1.0000&   -0.0180&   -1.0000&   -1.0000&    0.0198\\
   -0.5698 &  -0.0180&    1.0000&    0.0209&    0.0212&   -1.0000\\
   -0.1534 &  -1.0000&    0.0209&    1.0000&    1.0000&   -0.0227\\
   -0.1554 &  -1.0000&    0.0212&    1.0000&    1.0000&   -0.0229\\
    0.5705 &   0.0198&   -1.0000&   -0.0227&   -0.0229&    1.0000\\
\end{array}
\right].
\end{equation}
Suppose we start with a one-dimensional reduced basis corresponding to the first left singular vector of
$\statesnap$
\begin{equation*}
\nominalrightbasis=\rightbasiscoarse = \rightbasiscoarsei{1} = \left[
\begin{array}{c c c c c c}
  -0.2609   &-0.0348    &0.9390    &0.1240    &0.0463   &-0.1773
\end{array}
\right]^T.
\end{equation*}
Because the data nearly lie in a three-dimensional subspace of $\RR{6}$, the
optimal performance of a refinement scheme would yield small error after
splitting this one-dimensional basis into a basis of dimension three. Thus,
consider splitting $\rightbasiscoarse$ into three children using a tree that
follows the stated heuristic, i.e., is
characterized by $\children{1} = \{2,3,4\}$, $\elements{2} = \{1\}$,
$\elements{3} = \{3,6\}$, and $\elements{4} = \{2,4,5\}$. The resulting basis
becomes
\begin{equation*}
\rightbasisfine = \left[
\begin{array}{c c c c c c}
   -0.2609 &        0&         0&         0&         0&         0\\
         0 &        0&    0.9390&         0&         0&   -0.1773\\
         0 &  -0.0348&         0&    0.1240&    0.0463&         0
\end{array}
\right]^T.
\end{equation*}
The resulting projection error of the data is merely
$
\|\statesnap -
\rightbasisfine\left(\rightbasisfine\right)^+\statesnap\|_F/\|\statesnap\|_F =
0.0033.
$
By contrast, generating an alternative three-dimensional fine basis $\barrightbasisfine$ by splitting the basis using a (similar) tree characterized by
$\elements{2} = \{1\}$, $\elements{3} = \{3,5\}$, $\elements{3} = \{2,4,6\}$,
yields a much larger error of
$
\|\statesnap - \barrightbasisfine\left(\barrightbasisfine\right)^+\statesnap\|_F/\|\statesnap\|_F = 0.4948.
$

One way to identify these correlated variables is to employ $k$-means
clustering \cite{kmeans} after pre-processing the data by 1) normalizing observations of
each variable (to enable clustering to detect correlation), and 2) negating
the observation vector if the first observation is negative
(to enable clustering to detect anti-correlation). This is visualized in
Figure \ref{fig:clusteringProcessing} for the current example. Note that correlated and anti-correlated
variables have a small Euclidean distance between them after this processing;
this allows $k$-means clustering to identify them as a group.
\begin{figure}[htbp]
\centering
\subfigure[before processing (state variables labeled)]{
\includegraphics[width=.45\textwidth]{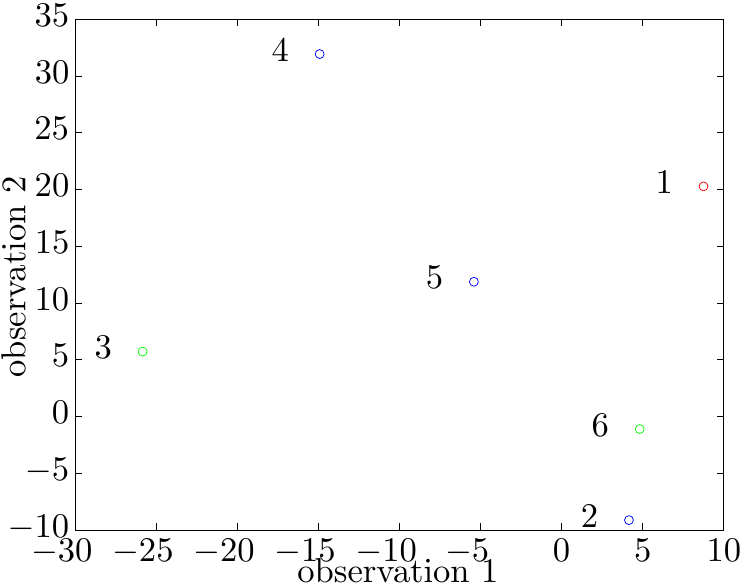}
}
\subfigure[after processing]{
\includegraphics[width=.45\textwidth]{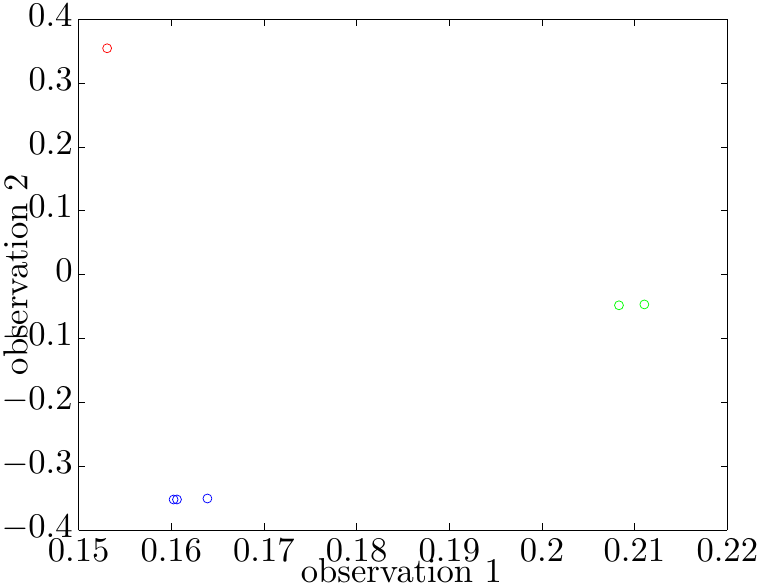}
\label{fig:snap1mesh378}
}
\caption{First two observations of the state variables (i.e., first two
columns of $\statesnap$) for the example in
Section \ref{sec:kmeans}. After processing these observations by normalization
and origin flipping, correlated and anti-correlated state variables are separated by
small geometric distances and can thus be grouped via clustering.}
\label{fig:clusteringProcessing}
\end{figure}
\end{example}
To this end, we construct the tree offline by recursively applying $k$-means
clustering to observations of the state variables (after reference subtraction,
normalization, and origin flipping). Algorithm \ref{alg:kmeans} describes the
method. The $\nobs$ observations of these variables are obtained from
snapshot data, which are often available, e.g.,  when the reduced basis is constructed
via proper orthogonal decomposition.

\begin{algorithm}[tbp]
\caption{Tree construction via recursive $k$-means clustering (offline)}
\begin{algorithmic}[1]\label{alg:kmeans}
\REQUIRE $\nobs$ snapshots of the reference-centered\footnotemark\ state in matrix form
$\statesnap\in\RR{\ndof\times \nobs}$, number of means $\kmeans$
\ENSURE child function $\childrenNo$, element function $\elementsNo$, and
number of nodes $\clustercount$
\FOR{$i=1,\ldots,\ndof$}
\STATE Normalize rows of $\statesnap$ to capture correlation by clustering
$\statesnaprowi{i}\leftarrow\statesnaprowi{i}/\|\statesnaprowi{i}\|$
\IF[Flip over origin to capture negative correlation by clustering]{$\statesnapij{i}{1} < 0$}
\STATE $\statesnaprowi{i}\leftarrow-\statesnaprowi{i}$
\ENDIF
\ENDFOR
\STATE Set root node to contain all elements $\elements{1} = \nat{\ndof}$.
\STATE Initialize recent-node set $\recentclusters \leftarrow \{1\}$ and node
count $\clustercount\leftarrow 1$.
\WHILE{$\card{\recentclusters} > 0$}
\STATE $\recentclustersold \leftarrow \recentclusters$, $\recentclusters \leftarrow \emptyset$

\FOR{$\icluster = 1,\ldots,\card{\recentclustersold}$}
\STATE Set splitting node to the $i$th element of the recent-node set $\clustertosplit
\leftarrow\recentclustersoldi{\icluster}$, where ordering is implied by $\geq$
on the natural numbers.
\IF[No elements to split]{ $\elements{\clustertosplit} = \emptyset$}
\STATE Continue
\ENDIF
\STATE Select snapshots of current elements $\statesnaptmpij{j}{k}
\leftarrow \statesnapij{\elementsNoij{\clustertosplit}{j}}{k}$, 
$j\in\nat{\card{\elements{\clustertosplit}}}$, $k\in\nat{\nobs}$
\STATE  $(\labelset_1,\ldots,\labelset_\ncluster)
=\kmeansmethod\left(\statesnaptmp,\kmeans\right)$, where
$\labelset_j\subset \nat{\card{\elements{\clustertosplit}}}$ denotes the set of elements in
cluster $j$, and $\ncluster$ denotes the number of non-empty clusters.
\IF[Cannot have only one child]{\ncluster = 1}
\FOR[Make all children into leaf nodes]{$j=1,\ldots,\card{\elements{\clustertosplit}}$}
\STATE $\labelset_j= j$
\ENDFOR
\ENDIF
\FOR{$\iChild = 1,\ldots,\ncluster$}
\STATE $\clustercount\leftarrow \clustercount + 1$
\STATE $\recentclusters\leftarrow \recentclusters \cup \clustercount$
\STATE $\elements{\clustercount} = \{\elementsNoij{\clustertosplit}{\iChild}\ |\
\iChild\in \labelset_{\iChild}\}$
\STATE $\childrenNoij{\clustertosplit}{\iChild} = \clustercount$
\ENDFOR
\ENDFOR
\ENDWHILE
\end{algorithmic}
\end{algorithm}
\footnotetext{This implies
that the reference state $\reference$ should be subtracted from the state
snapshots.}

\section{Dual-weighted residual error indicators}\label{sec:errorInd}
To compute error indicators for refinement, we propose a goal-oriented
dual-weighted residual methodology based on adjoint solves. It can be
considered a model-reduction adaptation of duality-based error-control methods
developed for differential equations
\cite{estep1995posteriori,pierce2000adjoint}, finite-element discretizations
\cite{babuvska1984post,becker1996weighted,rannacher1999dual,bangerth2003adaptive},
finite-volume discretizations
\cite{venditti2000adjoint,venditti2002grid,park2004adjoint}, and discontinuous
Galerkin discretizations \cite{lu2005posteriori,fidkowski2007simplex}.
Because the proposed method performs refinement online at the iteration level, it
requires error indicators associated with the error in ROM output at
iteration $\timestepit$, i.e., $\outputfun(\reference + \rightbasis\redstate^{{\timestepit}};\params)$. To
simplify notation in this section, we set $\reference = 0$ and write the
associated single solve
(Eq.~\eqref{eq:ROMdyn2} for a single iteration and parameter instance)
simply as 
\begin{equation}\label{eq:ROMsimple}
\rightbasis^T\modresidual(\rightbasis\redstate) =
0.
\end{equation}

 First, we approximate the output due to the (unknown) fine solution
$\redstatefine$ to first-order about the coarse solution $\redstatecoarse$:
\begin{equation} \label{eq:outputFirst}
\outputfun\left(
\rightbasisfine\redstatefine\right)
\approx
\outputfun\left(
\rightbasiscoarse\redstatecoarse\right) + \frac{\partial
\outputfun}{\partial \state}\left(
\rightbasiscoarse\redstatecoarse\right)\rightbasisfine\left(\redstatefine
- \coarsetofine\redstatecoarse\right),
\end{equation} 
where we have used Eq.~\eqref{eq:rightbasisfineSimple} to relate the coarse and fine bases.
 Similarly, we can approximate the fine residual to first order about the
 coarse solution as
 \begin{equation} 
 0=\rightbasisfineT\modresidual\left(
\rightbasisfine\redstatefine\right)\approx 
\rightbasisfineT\modresidual\left(
\rightbasiscoarse\redstatecoarse\right) + 
\rightbasisfineT\frac{\partial\modresidual}{\partial\state}\left(
\rightbasiscoarse\redstatecoarse\right)\rightbasisfine
\left(\redstatefine
- \coarsetofine\redstatecoarse\right).
\end{equation}
Solving for the state error yields
\begin{equation}\label{eq:resFirst}
\left(\redstatefine
- \coarsetofine\redstatecoarse\right)\approx
-
\left[\rightbasisfineT\frac{\partial\modresidual}{\partial\state}\left(
\rightbasiscoarse\redstatecoarse\right)\rightbasisfine\right]^{-1}
\rightbasisfineT\modresidual\left(
\rightbasiscoarse\redstatecoarse\right)
  \end{equation} 
Substituting \eqref{eq:resFirst} in \eqref{eq:outputFirst} yields
\begin{equation} \label{eq:outputFinal}
\outputfun\left(
\rightbasisfine\redstatefine\right)
-
\outputfun\left(
\rightbasiscoarse\redstatecoarse\right)
\approx
- \redadjointfineT\rightbasisfineT\modresidual\left(
\rightbasiscoarse\redstatecoarse\right).
\end{equation} 
where the fine adjoint solution $\redadjointfine\in\RR{\nrbfine}$ satisfies
\begin{equation}\label{eq:adjointFine}
\rightbasisfineT\frac{\partial\modresidual^{\timestepit}}{\partial\state}(
\rightbasiscoarse\redstatecoarse)^T\rightbasisfine\redadjointfine =
\rightbasisfineT\frac{\partial\outputfun}{\partial\state}\left(
\rightbasiscoarse\redstatecoarse\right)^T.
\end{equation}

Because we would like to avoid $\nrbfine$-dimensional solves associated with the fine basis
$\rightbasisfine$, we approximate $\redadjointfine$ as the prolongation of the
coarse adjoint solution
\begin{equation}\label{eq:fineadjoint}
\redadjointcoarsetofine=\coarsetofine\redadjointcoarse,
\end{equation}
where
$\redadjointcoarse$ satisfies
\begin{equation}\label{eq:adjointCoarse}
\rightbasiscoarseT\frac{\partial\modresidual^{\timestepit}}{\partial\state}(
\rightbasiscoarse\redstatecoarse)^T\rightbasiscoarse\redadjointcoarse =
\rightbasiscoarseT\frac{\partial\outputfun}{\partial\state}\left(
\rightbasiscoarse\redstatecoarse\right)^T
\end{equation}
Substituting the approximation $\redadjointcoarsetofine$ for $\redadjointfine$
in \eqref{eq:outputFinal} yields a cheaply computable error estimate
\begin{equation} \label{eq:outputFinalcheap}
\outputfun\left( \rightbasisfine\redstatefine\right) - \outputfun\left(
\rightbasiscoarse\redstatecoarse\right) \approx -
\redadjointcoarsetofineT\rightbasisfineT\modresidual\left(
\rightbasiscoarse\redstatecoarse\right).
\end{equation} 
The right-hand side can be bounded as
 \begin{equation} 
|\redadjointcoarsetofineT\rightbasisfineT\modresidual\left(
\rightbasiscoarse\redstatecoarse\right)|\leq\sum_{i\in\nat{\nrbfine}}\errorIndicatorFinei{i},
  \end{equation} 
	where the error indicators $\errorIndicatorFinei{i}\in\RR{}_+$,
	$i\in\nat{\nrbfine}$ are
\begin{equation}\label{eq:fineErrorEstimate}
\errorIndicatorFinei{i} =
|\redadjointcoarsetofinei{i}\rightbasisfineiT{i}\modresidual\left(
\rightbasiscoarse\redstatecoarse\right)|.
\end{equation}
Meyer and Matties
\cite{meyer2003efficient} also proposed a dual-weighted residual method for
reduced-order models. However, their approach was not applied to adaptive
refinement and did not consider a hierarchy of reduced bases; further, their
proposed dual solve was carried out on the full-order model, which is
infeasibly expensive for the present context.
\begin{remark}
Some mesh-refinement techniques \cite{venditti2000adjoint,venditti2002grid}
advocate computing refinement indicators that minimize the error in the
computable correction
$$\left(\redadjointfine-\redadjointcoarsetofine\right)^T\rightbasisfineT\modresidual\left(
\rightbasiscoarse\redstatecoarse\right).$$
To approximate this quantity, they employ prolongation operators of varying
fidelity, e.g., linear and quadratic interpolants. Such a strategy is not
straightforwardly applicable to the current context, as the prolongation operator
$\coarsetofine$ is exact.
\end{remark}

\begin{algorithm}[tbp]
\caption{Error estimates (online)}
\begin{algorithmic}[1]\label{alg:errorestimate}
\REQUIRE coarse reduced basis $\rightbasiscoarse$, coarse solution
$\redstatecoarse$
\ENSURE fine reduced basis $\rightbasisfine$, fine error-estimate vector $\errorIndicatorFine$
\STATE Solve coarse adjoint problem \eqref{eq:adjointCoarse} for
$\redadjointcoarse$.
\STATE Define prolongation operator $\coarsetofine$ via
Eq.~\eqref{eq:coarsetofine}.
\STATE Define fine reduced basis $\rightbasisfine$ via
Eq.~\eqref{eq:rightbasisfineSimple} and 
fine representation of adjoint solution $\redadjointcoarsetofine$ via Eq.~\eqref{eq:fineadjoint}
\STATE Compute fine error-estimate vector $\errorIndicatorFine$ via Eq.~\eqref{eq:fineErrorEstimate}
\end{algorithmic}
\end{algorithm}

\section{Adaptive $h$-refinement algorithm}\label{sec:adaptive}

We now return to the original objective of this paper: adaptively refine the
reduced-order model online. Algorithm \ref{alg:hrefine} describes our proposed
methodology for achieving this within a time-integration scheme. Step
\ref{step:solveRom} first computes the reduced-order-model solution satisfying
a tolerance $\romtol$. Then in Step \ref{step:fomNotConverged}, 
refinement occurs if the norm of the \emph{full-order residual} is above a
desired threshold $\fomtol$. Note that other (inexpensive) error indicators could be used to
flag refinement, e.g., error surrogates \cite{drohmannRomes}. Refinement continues until this full-order
tolerance is satisfied; note that any tolerance can be reached, as a
completely split basis yields a reduced-order model equivalent to the
full-order model (see Section \ref{sec:tree}).  Finally, Step
\ref{step:resetbasis} resets the basis every
$\resetFreq$ time iterations. This ensures 1) the basis does not grow monotonically,
and 2) work performed to refine the basis can be amortized over subsequent
time
steps, where the solution is unlikely to significantly change.  Note that if
Step \ref{step:solveRom} entails an iterative solve (e.g., Newton), then the
pre-refinement solution can be employed as an initial guess.

\begin{algorithm}[tbp]
\caption{Adaptive $h$-refinement (online)}
\begin{algorithmic}[1]\label{alg:hrefine}
\REQUIRE iteration ${\timestepit}$, basis $\rightbasis$, ROM solver tolerance $\romtol$, FOM solver tolerance $\fomtol$
\ENSURE updated basis $\rightbasis$, generalized state $\redstate^{\timestepit}$
\STATE\label{step:solveRom} Compute ROM solution $\redstate^{\timestepit}$ satisfying
$\|\rightbasis^T\modresidual^{\timestepit}(\reference +
\rightbasis\redstate^{{\timestepit}};\params)\|\leq \romtol$.
\IF{FOM not converged $\|\modresidual^{\timestepit}(\reference +
\rightbasis\redstate^{{\timestepit}};\params)\|>
\fomtol$}\label{step:fomNotConverged}
\STATE Refine basis  via
Algorithm \ref{alg:refine}:
$\rightbasis\leftarrow\mathrm{Refine}\left(\rightbasis,\redstate^{\timestepit}\right)$.
\STATE Return to Step \ref{step:solveRom}.
\ENDIF
\IF{$\mod({\timestepit},\resetFreq) =0$}
\STATE\label{step:resetbasis} Reset basis $\rightbasis \leftarrow \nominalrightbasis$.
\ENDIF
\end{algorithmic}
\end{algorithm}

Algorithm \ref{alg:refine} describes the proposed method for refining the
basis using the refinement mechanism and error indicators presented in
Sections \ref{sec:refinement} and \ref{sec:errorInd}, respectively. 
Appendix \ref{app:multipleTrees} describes a more sophisticated approach
wherein the basis vectors are not split into all possible children; the
children are separated into groups, each of which contributes roughly
the same 
fraction of that vector's error.

First, Step \ref{step:dwrerror} of Algorithm \ref{alg:refine} computes error estimates for the fine basis
(i.e., current basis with all vectors split into all possible children) using
the dual-weighted residual approach. Step \ref{step:flagvectors} marks
the parent basis vectors to refine: those with above-average error
contribution from its children. Steps \ref{step:beginChildSplit}--\ref{step:endChildSplit} split the parent vector
$i$ into vectors corresponding to its $\nchildren{i}$ children according 
to the defined tree.
Steps \ref{step:vectoraddstart}--\ref{step:vectoraddend} update the reduced basis and tree nodes.
Because this split does not guarantee a full-ranks basis, Step \ref{step:qr}
performs an efficient QR factorization with column pivoting to identify
`redundant' basis vectors. Step \ref{step:fullrankstart} subsequently removes
these vectors from the basis and Step \ref{step:fullrankend} performs the
necessary bookkeeping for the tree nodes.

\begin{algorithm}[tbp]
\caption{Refine (online)}
\begin{algorithmic}[1]\label{alg:refine}
\REQUIRE initial basis $\rightbasis$, reduced solution $\redstate$
\ENSURE refined basis $\rightbasis$
\STATE\label{step:dwrerror} Compute fine error-estimate vector  and fine reduced
basis  via Algorithm \ref{alg:errorestimate}:\\
$\left(\errorIndicatorFine,\rightbasisfine\right)\leftarrow\mathrm{Error\
estimates}\left(\rightbasis,\redstate\right)$.
\STATE Put local error estimates in parent--child format
$\locerror{i}{j} = \errorIndicatorFinei{\vectorparenttochild{i}{j}}$,
$i\in\nat{\nrb}$, $j\in\nat{\nchildren{i}}$.
\STATE\label{step:flagvectors} Mark basis vectors to refine $\flaggedvectors= \{i\ |\
\sum_j\locerror{i}{j} \geq 1/p\sum_{kj}\locerror{k}{j}\}$
\FOR[Split $\rightbasisi{i}$ into $\nchildren{i}$ vectors]{$i\in\flaggedvectors$}
\FOR{$k\in\nat{\nchildren{i}}$}\label{step:beginChildSplit}
\STATE $\extrabasis{k} =\rightbasisfinei{\vectorparenttochild{i}{k}}$
\STATE $\tmpnode_k = \childrenNoij{\node_i}{k}$
\ENDFOR\label{step:endChildSplit}
\STATE\label{step:vectoraddstart} $\rightbasisi{i}\leftarrow\extrabasis{1}$,
$\node_{i}\leftarrow\tmpnode_{1}$
\FOR{$k= 2,\ldots\nchildren{i}$}
\STATE $\rightbasisi{\nrb +
k-1}\leftarrow\extrabasis{k}$,
$\node_{\nrb +
k-1}\leftarrow\tmpnode_{k}$,
\ENDFOR\label{step:vectoraddend}
\ENDFOR
\STATE\label{step:qr} Compute thin QR factorization with column pivoting $ \rightbasis = \Q\R$, 
$\R\Perm =\Qs\Rs$.
\STATE\label{step:fullrankstart} Ensure full-rank matrix $ \rightbasis \leftarrow\rightbasis
\left[\Permi{1}\ \cdots\
\Permi{\rankbasis}\right]$, where $\rankbasis$ denotes the numerical rank of $\R$. 
\STATE \label{step:fullrankend}Update tree 
$\left[\node_{1}\ \cdots\
\node_{\rankbasis}\right]\leftarrow\left[\node_{1}\ \cdots\ \node_{\nrb}\right]\left[\Permi{1}\ \cdots\
\Permi{\rankbasis}\right]$.
\end{algorithmic}
\end{algorithm}

\section{Numerical experiments: parameterized inviscid Burgers'
equation}\label{sec:experiments}

We assess the method's performance on the parameterized inviscid Burgers'
equation. While simple, this problem is particularly challenging for reduced-order
models. This arises from the fact that ROMs approximate the solution as a
linear combination of spatially fixed reduced-basis functions; as such, they work well
when the dynamics are primarily Eulerian, i.e., are fixed with respect to the
underlying grid. However, when the dynamics are Lagrangian in nature and exhibit motion
with respect to the underlying grid (e.g., moving shocks), reduced-order
models generally fail to capture the critical phenomenon at every time step and
parameter instance.

We employ the problem setup described in Ref.~\cite{tpwl}.
Consider the parameterized initial boundary value problem
\begin{align}
\label{eq:burgers}
\frac{\partial u(x,\timeVariable)}{\partial \timeVariable} + \frac{1}{2}\frac{\partial
\left(u^2\left(x,\timeVariable\right)\right)}{\partial x} &= 0.02e^{\param{2}x}\\
u(0,\timeVariable) &= \param{1}, \ \forall \timeVariable>0\\
\label{eq:burgersLast}u(x,0) &= 1, \ \forall x\in\left[0,~100\right],
\end{align}
where $\param{1}$ and $\param{2}$ are two real-valued input variables. Godunov's scheme
discretizes the problem, which leads to a finite-volume
formulation consistent with the original formulation in Eq.~\eqref{eq:FOMdyn}. The one-dimensional domain is discretized using a grid with
$\burgersSizep$
nodes corresponding to coordinates coordinates $x_i = i \times
(100/\burgersSize)$,
$i=0,\ldots, \burgersSize$.  Hence, the resulting full-order model is of dimension
$\ndof=\burgersSize$.
The solution $u(x,\timeVariable)$ is computed in the time interval
$\timeVariable\in\left[0,50\right]$ using a uniform computational time-step size
$\Delta t = 0.05$, leading to $\nt = 1000$ total time steps.

For simplicity, we employ a POD--Galerkin ROM. During the
offline stage,
snapshots of the state are collected for the first $\ntTrain$ time steps at
training inputs. Then, the initial condition is subtracted from these
snapshots, and they are concatenated column-wise to generate the snapshot
matrix.
Finally, the thin singular value decomposition of the snapshot matrix is
computed, and the initial reduced basis $\nominalrightbasis$ is set to the
first $\nrbinitial$ left singular
vectors. During the online stage, a Galerkin projection is employed using this
reduced basis.
For all experiments, the initial condition is
set to the reference condition, i.e., $\reference = \state^0$.
For $h$-adaptivity, we set the number of means to $\kmeans=10$
in Algorithm \ref{alg:kmeans}.
For Algorithm \ref{alg:errorestimate}, the output of interest is set to
the residual norm, i.e., 
$\outputfun(\state^{\timestepit};\params) = \|
\residual^{\timestepit}(\state^{{\timestepit}};\params)\|_2^2$. 
For Algorithm \ref{alg:hrefine}, the ROM tolerance is set to $\romtol =
5\times 10^{-3}$.\footnote{For the ROMs without
adaptivity, the ROM convergence tolerance is set to $\romtol = 1\times
10^{-5}$.} The basis-reset
frequency $\resetFreq$ will vary during the experiments. Step \ref{step:solveRom}
incurs a Newton solve; when refinement has occurred, the initial guess is
set to the converged solution from the
previous refinement level. Finally, the
experiments
employ the (more complex) Refine method defined by
Algorithm \ref{alg:refinemult} with a child-partition factor
$\npartition=2$.

Note that because the residual operator is nonlinear in the state, a
projection alone is insufficient to generate computational savings over the
full-order model. Future work will address extending the proposed
$h$-refinement method to ROMs equipped with a complexity reduction mechanism
such as empirical interpolation or gappy POD.

\subsection{Fixed inputs}
For this example, the input parameters are set to $\param{1}=3$ and
$\param{2}=0.02$. However, the problem can be considered to be predictive, as
we only collect snapshots in the time interval $\tTrain\in\left[0,7.5\right]$,
i.e., for the first $\ntTrain = 150$ time steps. This choice is made to
introduce a significant challenge for the ROM: while the (unrefined) reduced
basis captures discontinuities that arise in the first 150 time
steps, it will not capture such discontinuities that arise outside of this
time interval.\footnote{Note that the refinement method can also be
applied when the original reduced basis captures all relevant online
phenomena; however, the need for \emph{a posteriori} refinement is weaker
in this case.}

\begin{figure}[htbp]
\centering
\subfigure[no adaptivity, no truncation ($\nrb=150$)]{
\includegraphics[width=6.5cm]{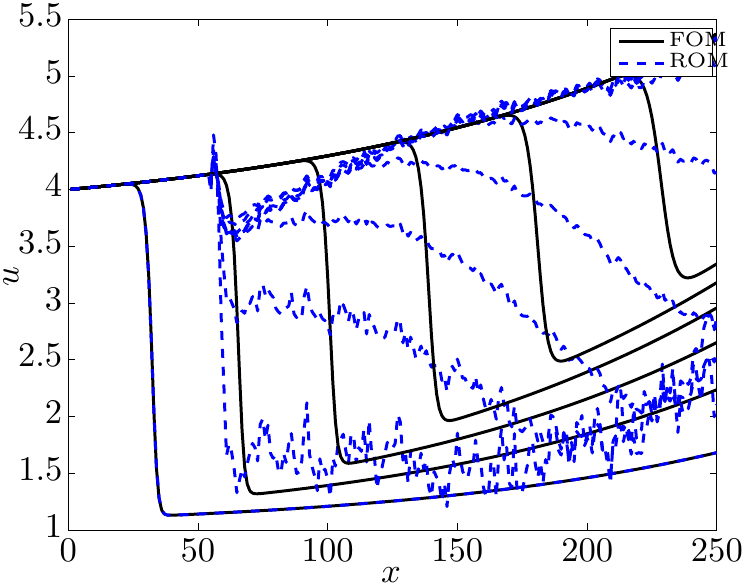} }
\subfigure[$h$-adaptivity, $\nrbinitial = 10$, $\resetFreq=50$]{
\includegraphics[width=6.5cm]{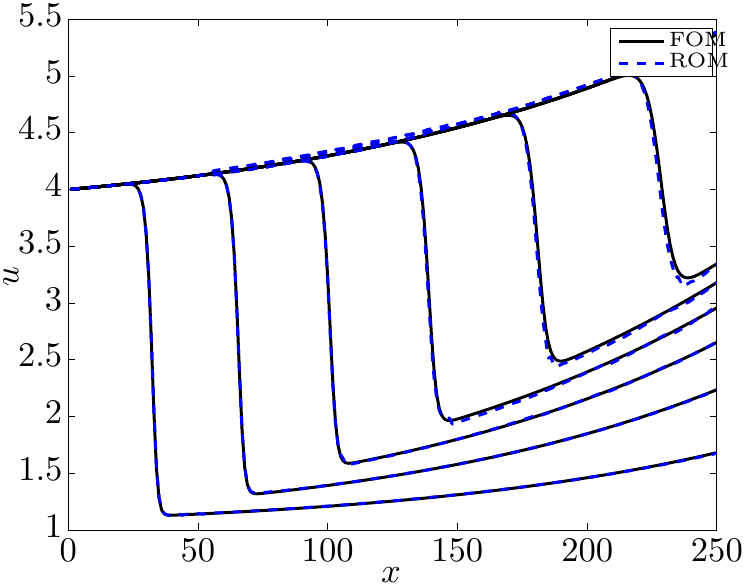} }
\caption{Comparison of solutions computed by POD--Galerkin with and without
$h$-adaptivity for the fixed-inputs case.}
\label{fig:compareSolutions}
\end{figure}

Table \ref{tab:resultsNominal} reports results for typical POD--Galerkin ROMs of
differing dimensions, as well as results for the proposed $h$-refinement
method with different parameters and a FOM tolerance in Algorithm
\ref{alg:hrefine} of $\fomtol = 0.05$. Here, the relative error is defined as $$
\mathrm{relative\ error}=\frac{1}{\nt}\sum_{\timestepit=1}^\nt\|u_{\mathrm{FOM}}(\cdot,\timeVariable^\timestepit)-\urom(\cdot,\timeVariable^\timestepit)\|_{L_2}/
\|u_{\mathrm{FOM}}(\cdot,\timeVariable^\timestepit)\|_{L_2}.  $$ Figure
\ref{fig:compareSolutions} compares the solutions predicted by POD--Galerkin
with no basis truncation (i.e., $\nrb=150$) and that of the proposed method
with an initial basis size of $\nrbinitial = 10$ with
$\nominalrightbasis\in\RR{\ndof\times\nrbinitial}$ and a basis-reset
frequency of $\resetFreq=50$. 

\begin{table}[htd] 
\small
\centering
\begin{tabular}{|c||c|c|c||c|c|c|c|c|}
\hline
&\multicolumn{3}{|c||}{no adaptivity} & \multicolumn{5}{|c|}{$h$-adaptivity}\\
\hline
initial basis dimension $\nrbinitial$ & 10 & 45 & 150 & 5 &10 &  20 & 10 & 10\\
basis-reset frequency $\resetFreq$& & & & 50 & 50 &
50 &
100 & 25\\
\hline
average basis dimension & \multirow{2}{*}{10} & \multirow{2}{*}{45} &
\multirow{2}{*}{150} & \multirow{2}{*}{41.4}& \multirow{2}{*}{44.3}  &
\multirow{2}{*}{58} & \multirow{2}{*}{73} & \multirow{2}{*}{37}\\
per Newton iteration $\nrbavg$& & & & &   & & & \\
average number of Refine & & & &  \multirow{2}{*}{$0.20$} &
\multirow{2}{*}{$0.19$} &\multirow{2}{*}{$0.14$} & \multirow{2}{*}{$0.13$} &
\multirow{2}{*}{$0.28$}\\
calls per time step & & & &  & && & \\
\hline
relative error ($\%$) & $45.8$ & $43.9$ & $8.5$ & $0.3$ & $0.5$ & $0.2$ &
$0.2$ &
$0.3$\\
online time (seconds)& $1.4$ & $2.14$ & $5.77$ & $5.53$ &$4.63$ &
$7.27$ & $6.90$ &
$7.46$\\
\hline
\end{tabular}
\caption{Comparison between POD--Galerkin ROMs without refinement and with
$h$-adaptive refinement for the fixed-inputs case.} 
\label{tab:resultsNominal} 
\end{table}

First, note that the reduced-order model is highly inaccurate (even when the
basis is not truncated) unless equipped with $h$-adaptivity. The reason for
this is simple: the training has not captured the flow regime with shock
locations past approximately $x=60$. This illustrates a powerful capability of
the proposed $h$-adaptation methodology: it enables ROMs to be incrementally
refined to capture previously unobserved phenomena. In fact, the average basis
dimension (per Newton iteration) for the best-performing $h$-adaptive ROM ($\nrbinitial=10$,
$\resetFreq=50$) is only $\nrbavg=44.3$, which is smaller than the basis dimensions for
ROMs without adaptivity ($\nrb=45$ and $\nrb=150$) that yield much higher errors
($43.9\%$ and $8.5\%$, respectively).

Second, adaptation parameters $\nrbinitial$ and $\resetFreq$ both lead to a
performance tradeoff. When $\nrbinitial$ is small, it leads to smaller average
basis sizes $\nrbavg$. However, it increases the number of Refine calls per
time step, as the smaller basis must be refined more times to achieve desired
accuracy. Similarly, resetting the basis more frequently (smaller
$\resetFreq$) leads to a smaller $\nrbavg$, but more average refinement steps.
As such, an intermediate value of both parameters leads to the shortest online
evaluation time.

Finally, notice that the online evaluation time for the adaptive ROM with an
average basis size of $\nrbavg=44.3$ is roughly twice that of a non-adaptive ROM
with roughly the same basis size $\nrb=45$. This discrepancy in evaluation time
can be attributed to the overhead in performing the adaptation. For larger
problem sizes, one would expect this overhead to be smaller relative to the
total online evaluation time.

\begin{figure}[htbp]
\centering
\subfigure[$\fomtol = 0.35$]{
\includegraphics[width=6.5cm]{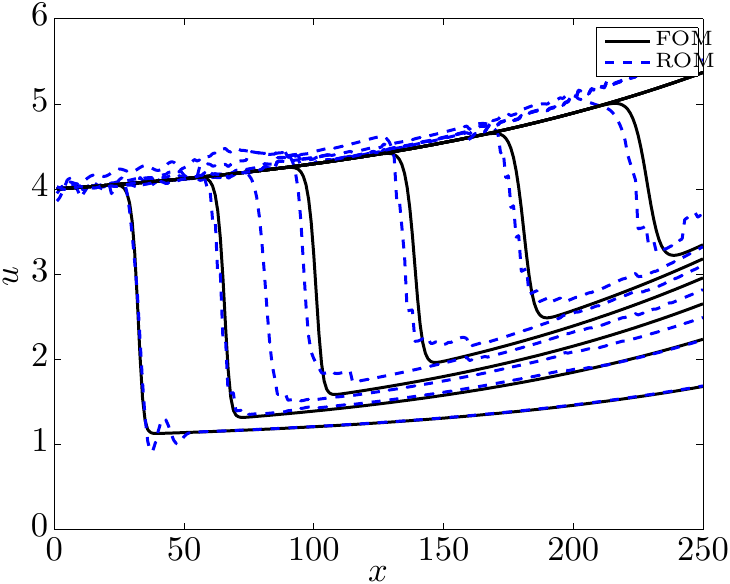} }
\subfigure[$\fomtol = 0.05$]{
\includegraphics[width=6.5cm]{kevin1000C} }
\subfigure[$\fomtol = 0.01$]{
\includegraphics[width=6.5cm]{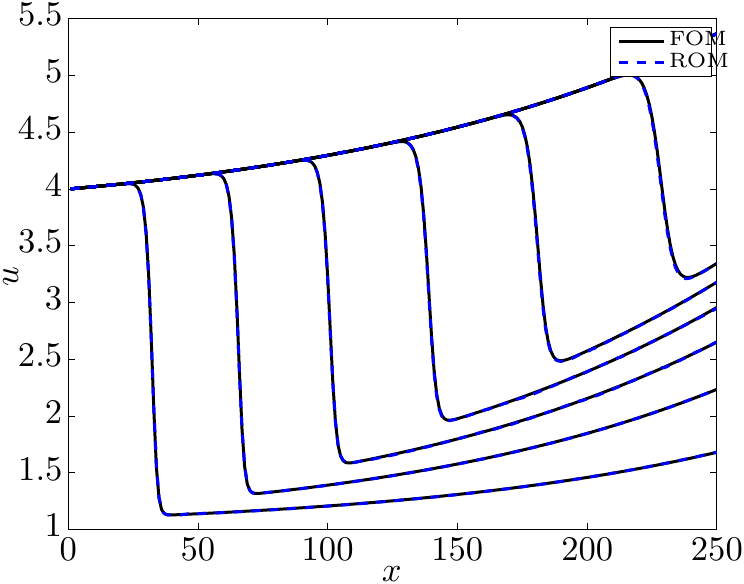} }
\caption{Comparison of solutions computed by $h$-adaptive POD--Galerkin for
different full-order-model tolerances $\fomtol$ for the fixed-inputs case.}
\label{fig:fomtol}
\end{figure}

Next, we assess the performance of the $h$-refinement method as the
full-order-model
tolerance $\fomtol$ in Algorithm \ref{alg:hrefine} varies. Table \ref{tab:fomtol} and Figure \ref{fig:fomtol}
report the results. As expected, the proposed method allows the ROM to achieve
any of the prescribed tolerances. As the tolerance becomes more rigorous, the
ROM solution improves; however, it does so at increased computational cost, as
both the average basis dimension $\nrbavg$ and number of Refine calls per time
step increase to satisfy the requirement. 

\begin{table}[htd] 
\small
\centering
\begin{tabular}{|c||c|c|c|}
\hline
&$\fomtol = 0.35$ & $\fomtol = 0.05$ & $\fomtol = 0.01$\\
\hline
average basis dimension & \multirow{2}{*}{33.6} & \multirow{2}{*}{44.2507} &
\multirow{2}{*}{53.9}\\
per Newton iteration $\nrbavg$& & &  \\
average number of Refine & \multirow{2}{*}{0.115} & \multirow{2}{*}{0.189} &
\multirow{2}{*}{0.212} \\
calls per time step & & &  \\
\hline
relative error ($\%$) & $12.2$ & $0.51$ & $0.078$\\
online time (seconds)& $4.61$ & $4.63$ & $7.64$\\
\hline
\end{tabular}
\caption{Effect of full-order-model tolerance $\fomtol$ on 
$h$-adaptive refinement for $\nrbinitial=10$ and $\resetFreq=50$ for the
fixed-inputs case.} 
\label{tab:fomtol} 
\end{table}

\subsection{Input variation}
For this experiment, we assess the proposed methodology in an input-varying
scenario. In particular, the offline stage collects snapshots in the time
interval $\tTrain\in\left[0,2.5\right]$ for the training set
$\{\params^1,\ldots,\params^3\}$ described in 
Table~\ref{tab:burgersGlobalInputs}, which is constructed by uniformly
sampling the input space along $(\param{1},\param{2}) = (3\alpha, 0.02\alpha)$,
$\alpha\in\left[1,3\right]$. 
\begin{table}[tb] 
	\caption{Offline and online inputs for the inviscid Burgers equation}
  \label{tab:burgersGlobalInputs} 
 \centering 
 \begin{tabular}{|c||c|c|c|c|} 
  \hline 
\multirow{2}{*}{Input variables} & Training point& Training point& Training
point& Online point\\
                                    &$\params^1$        &$\params^2$
																		&$\params^3$        & $\params^\star$\\
  \hline 
$\param{1}$ & 3 & 6 &9 & 4.5 \\
$\param{2}$ & 0.02& 0.05 &0.075 &  0.038\\
  \hline 
  \end{tabular} 
  \end{table} 

\begin{figure}[htbp]
\centering
\subfigure[no adaptivity, no truncation ($\nrb=150$)]{
\includegraphics[width=6.5cm]{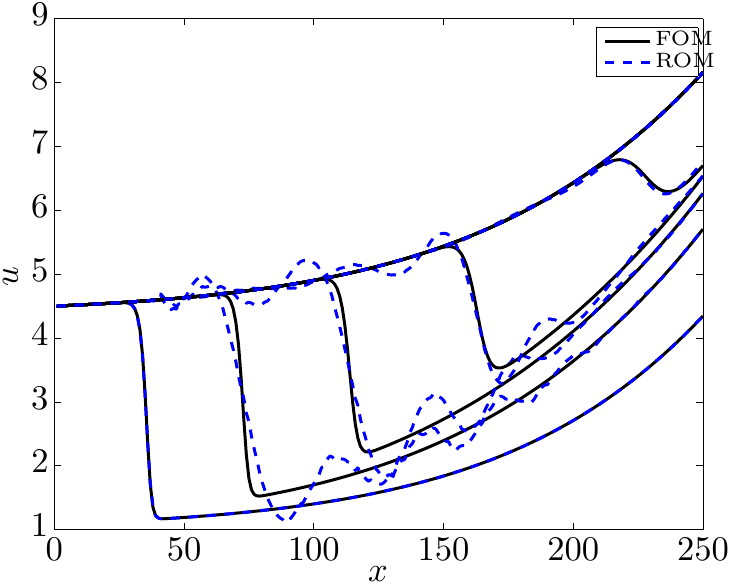} }
\subfigure[$h$-adaptivity, no truncation ($\nrbinitial = 20$), $\resetFreq=100$]{
\includegraphics[width=6.5cm]{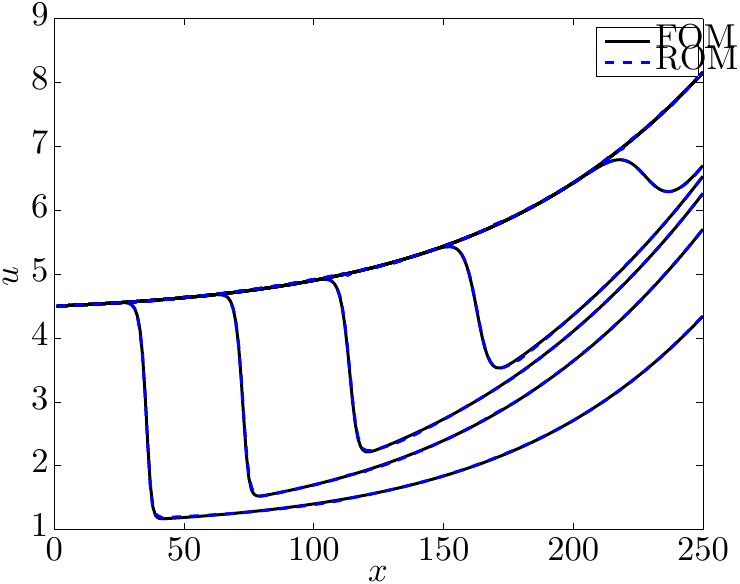} }
\caption{Comparison of solutions computed by POD--Galerkin with and without
adaptivity for the varying-inputs case.}
\label{fig:compareSolutionsParamVarying}
\end{figure}

Figure \ref{fig:compareSolutionsParamVarying} and Table
\ref{tab:resultsPredictive} report the results for this experiment. The same phenomena are prevalent as were
apparent in the previous experiment. The primary difference is that the
POD--Galerkin model without adaptivity performs better than previously (due to
more informative snapshots). However,
$h$-adaptivity is still required to drive errors below $1\%$. Note that the
proposed method compensated for an unsophisticated uniform-sampling of the
input space. The method would still be applicable for more rigorous (e.g.,
POD--Greedy \cite{haasdonk2008FV}) sampling methods, which would lead to a more robust initial basis
$\nominalrightbasis$ and reduce the burden of $h$-adaptivity to generate
accurate results.

\begin{table}[htd] 
\small
\centering
\begin{tabular}{|c||c|c|c||c|c|c|c|c|}
\hline
&\multicolumn{3}{|c||}{no adaptivity} & \multicolumn{5}{|c|}{$h$-adaptivity}\\
\hline
initial basis dimension $\nrbinitial$ & 10 & 78 & 150 & 5 &20 &  30 & 20 & 20\\
basis-reset frequency $\resetFreq$& & & & 100 & 100 &
100 &
200 & 50\\
\hline
average basis dimension & \multirow{2}{*}{10} & \multirow{2}{*}{78} &
\multirow{2}{*}{150} & \multirow{2}{*}{69.8}& \multirow{2}{*}{77.2}  &
\multirow{2}{*}{87.6} & \multirow{2}{*}{130.6} & \multirow{2}{*}{65.6}\\
per Newton iteration $\nrbavg$& & & & &   & & & \\
average number of Refine & & & &  \multirow{2}{*}{$0.20$} &
\multirow{2}{*}{$0.072$} &\multirow{2}{*}{$0.07$} & \multirow{2}{*}{$0.044$} &
\multirow{2}{*}{$0.11$}\\
calls per time step & & & &  & && & \\
\hline
relative error ($\%$) & $41.8$ & $1.7$ & $1.4$ & $0.22$ & $0.14$ & $0.45$ &
$0.53$ &
$0.70$\\
online time (seconds)& $1.75$ & $3.54$ & $8.55$ & $6.41$ &$6.06$ &
$8.11$ & $9.11$ &
$8.78$\\
\hline
\end{tabular}
\caption{Comparison between POD--Galerkin ROMs without refinement and with
$h$-adaptive refinement for the input-variation case.} 
\label{tab:resultsPredictive} 
\end{table}

\section{Conclusions}

This work has presented an adaptive $h$-refinement method for reduced-order
models. Key components include 1) an $h$-refinement mechanism based on basis
splitting and tree structure constructed via $k$-means clustering, 2)
dual-weighted residual error indicators, and 3) an adaptive algorithm to
moderate when and how to perform the refinement. In contrast to existing
\emph{a priori} adaptive methods, the proposed technique provides a mechanism
to improve the ROM solution \emph{a posteriori}. As opposed to existing
\emph{a posteriori} methods, the proposal does so without incurring
any large-scale operations. Numerical examples on the inviscid Burgers equation
highlighted the method's ability to accurately predict phenomena not present
in the training data used to construct the reduced basis.

Future research directions include incorporating complexity reduction (e.g.,
empirical interpolation, gappy POD) into the
refinement process. In particular, as the reduced basis is refined, 
sample points (and dual reduced-basis vectors) must be added in a systematic
way to ensure the reduced-order model remains solvable.
Similar to the manner
in which the tree defining the (complete) splitting mechanism is constructed
offline, one could generate a hierarchy of these sample points offline from the
training data, e.g., by executing \cite[Algorithm 3]{carlbergJCP}
for $n_s$ equal to the number of nodes in the mesh.
In addition, it would be interesting to incorporate a more sophisticated
adaptive \emph{coarsening} technique (compared to the simple basis-resetting mechanism in
Step \ref{step:resetbasis} of Algorithm \ref{alg:hrefine}); for example, one could combine
basis vectors whose generalized coordinates are strongly correlated (or
anti-correlated) over recent time steps. Further, it would be interesting to
pursue adaptive $p$-refinement methods, wherein other basis vectors (e.g.,
truncated POD vectors, discrete wavelets) with possibly global support are added from a library to enrich
the reduced basis. In addition, it would be useful to pursue alternative
tree-construction methods that satisfy Conditions 
\ref{req:globalSupport}--\ref{req:fullSplit} of Section \ref{sec:tree}.
Assessing the effect of the proposed refinement method on ROM stability
would also constitute an interesting investigation.
Finally, it would be advantageous to incorporate
Richardson extrapolation in the refinement method to better approximate the
outputs of interest; however, this requires knowledge of the convergence rate
of the reduced-order model with respect to adding basis vectors.

\appendix
\section{Refinement algorithm with multiple trees}\label{app:multipleTrees}
This section presents a more sophisticated refinement mechanism than that that
presented in Section \ref{sec:adaptive}. In particular, when a vector is
flagged for refinement, it is not necessarily split into all its children.
Rather, its children are separated into groups, each of which contributes roughly
the same 
fraction $\npartition$ of the total error for that parent vector. This avoids
over-refinement when the number of children is relatively large. However, this
leads to an increase in required bookkeeping, as the tree structure changes
when children merge:
the tree
must be altered and separately maintained for each vector.
Thus, each basis vector $\rightbasisi{i}$, $i=1,\ldots,\nrb$  will be
characterized by its own tree $\childrenNo_i$, $\elementsNo_i$ with $\nnode_i$
nodes, as well as a
node on that tree $\node_i\in\nat{\nnode_i}$. 

Algorithm \ref{alg:refinemult} describes the modifications needed to Algorithm
\ref{alg:refine} to enable this feature. Key modifications include the
following. Steps
\ref{step:whilebegin2}--\ref{step:whileend2}
separate the children of the parent vector's tree node $\node_i$ into groups;
the resulting maintenance of the tree structures is
performed in Steps \ref{step:treealterchildren2}--\ref{step:treealterelements2}.\footnote{Only the
lower levels of the tree must be updated, as the current methodology never
traverses up a tree.} In steps
\ref{step:vectoraddstart2}--\ref{step:vectoraddend2},
not only is the basis updated, but the trees are as well.
Finally, Step \ref{step:fullrankend2} performs the
necessary bookkeeping for the tree structures due to the removal of redundant
basis vectors.

\begin{algorithm}[tbp]
\caption{Refine (child grouping) (online)}
\begin{algorithmic}[1]\label{alg:refinemult}
\REQUIRE initial basis $\rightbasis$, reduced solution $\redstate$,
child-partition factor $\npartition\leq 1$
\ENSURE refined basis $\rightbasis$
\STATE\label{step:dwrerror2} Compute fine error-estimate vector  and fine reduced
basis  via Algorithm \ref{alg:errorestimate}:\\
$\left(\errorIndicatorFine,\rightbasisfine\right)\leftarrow\mathrm{Error\
estimates}\left(\rightbasis,\redstate\right)$.
\STATE Put local error estimates in parent--child format
$\locerror{i}{j} = \errorIndicatorFinei{\vectorparenttochild{i}{j}}$,
$i\in\nat{\nrb}$, $j\in\nat{\nchildren{i}}$.
\STATE\label{step:flagvectors2} Mark basis vectors to refine $\flaggedvectors= \{i\ |\
\sum_j\locerror{i}{j} \geq 1/p\sum_{kj}\locerror{k}{j}\}$
\FOR[Split $\rightbasisi{i}$ into $k$ vectors]{$i\in\flaggedvectors$}
\STATE $\nrb\leftarrow \mathrm{dim}\left(\range{\rightbasis}\right)$
\STATE \label{step:forbegin}Initialize additional-vector count $k\leftarrow 0$ and handled child-node set  $\childgroupNo \leftarrow\emptyset$
\WHILE[Divide child nodes into groups with roughly equal
error]{$\childgroupNo\not = \nat{\nchildren{i}}$}\label{step:whilebegin2}
\STATE $\childgroup{k} = \arg\min_{\childgroupvar\subset
\childgroupset}\mathrm{card}\left(\childgroupvar\right)$, where
$\childgroupset = \{\childgroupvar\subset \nat{\nchildren{i}}
\setminus\childgroupNo\ |\
\sum_{j\in z}\locerror{i}{j}\geq \npartition\sum_{j}\locerror{i}{j} \}$.

\IF{$\childgroup{k} = \emptyset$}
\STATE Take all remaining children $\childgroup{k} =\nat{\nchildren{i}}\setminus\childgroupNo $
\ENDIF
\STATE $\extrabasis{k} =
\sum_{j\in \childgroup{k}}\rightbasisfinei{\vectorparenttochild{i}{j}}$
\STATE Update tree: $\tmpchildreniNo{k}\leftarrow\childreniNo{i}$,
$\tmpelementsiNo{k}\leftarrow\elementsiNo{i}$
\IF[Use the same tree]{$\mathrm{card}(\childgroup{k})=1$}
\STATE $\tmpnode_k = \childrenij{i}{\node_i}{\childgroup{k}}$
\ELSE[Alter the tree]
\STATE $\tmpnode_k =\node_i$
\STATE 
\label{step:treealterchildren2}$\tmpchildreni{k}{\tmpnode_k} = \{\childrenij{i}{\node_i}{k}\ |\
k\in\childgroup{k}\} $
\STATE 
\label{step:treealterelements2}$\tmpelementsi{k}{\tmpnode_k} = \bigcup_{k\in\tmpchildreni{k}{\node_i}}\elementsi{i}{k}$
\ENDIF
\STATE $k\leftarrow k+1$,
$\childgroupNo\leftarrow\childgroupNo\cup\childgroup{k}$
\ENDWHILE\label{step:whileend2}
\STATE\label{step:vectoraddstart2} $\rightbasisi{i}\leftarrow\extrabasis{0}$,
$\childreniNo{i}\leftarrow\tmpchildreniNo{0}$, 
$\elementsiNo{i}\leftarrow\tmpelementsiNo{0}$
$\node_{i}\leftarrow\tmpnode_{0}$
\FOR{$\extravecCounter = 1,\ldots,k$}
\STATE $\rightbasisi{\nrb +
\extravecCounter}\leftarrow\extrabasis{\extravecCounter}$,
$\childreniNo{\nrb +
\extravecCounter}\leftarrow\tmpchildreniNo{\extravecCounter}$, 
$\elementsiNo{\nrb +
\extravecCounter}\leftarrow\tmpelementsiNo{\extravecCounter}$,
$\node_{\nrb +
\extravecCounter}\leftarrow\tmpnode_{\extravecCounter}$
\ENDFOR\label{step:vectoraddend2}
\ENDFOR
\STATE\label{step:qr2} Compute thin QR factorization with column pivoting $ \rightbasis = \Q\R$, 
$\R\Perm =\Qs\Rs$.
\STATE\label{step:fullrankstart2} Ensure full-rank matrix $ \rightbasis \leftarrow\rightbasis
\left[\Permi{1}\ \cdots\
\Permi{\rankbasis}\right]$, where $\rankbasis$ denotes the numerical rank of $\R$. 
\STATE \label{step:fullrankend2}Update tree 
$\left[\childreniNo{1}\ \cdots\
\childreniNo{\rankbasis}\right]\leftarrow\left[\childreniNo{1}\ \cdots\ \childreniNo{\nrb}\right]\left[\Permi{1}\ \cdots\
\Permi{\rankbasis}\right]$;\\
$\left[\elementsiNo{1}\ \cdots\
\elementsiNo{\rankbasis}\right]\leftarrow\left[\elementsiNo{1}\ \cdots\ \elementsiNo{\nrb}\right]\left[\Permi{1}\ \cdots\
\Permi{\rankbasis}\right]$;\\
$\left[\node_{1}\ \cdots\
\node_{\rankbasis}\right]\leftarrow\left[\node_{1}\ \cdots\ \node_{\nrb}\right]\left[\Permi{1}\ \cdots\
\Permi{\rankbasis}\right]$.
%
%
\end{algorithmic}
\end{algorithm}

\section*{Acknowledgments}
The author acknowledges Matthew Zahr for providing the model-reduction testbed
that was modified to generate the numerical results, Seshadhri Comandur for
helpful discussions related to tree construction and clustering, and the
anonymous reviewers for providing insightful remarks and suggestions.  This research was
supported in part by an appointment to the Sandia National Laboratories Truman
Fellowship in National Security Science and Engineering, sponsored by Sandia
Corporation (a wholly owned subsidiary of Lockheed Martin Corporation) as
Operator of Sandia National Laboratories under its U.S.  Department of Energy
Contract No. DE-AC04-94AL85000. 
\bibliography{references}
\bibliographystyle{aiaa}
\end{document}